\documentclass{amsart}
\usepackage{eurosym}
\usepackage{amssymb,amsmath,amsfonts}
\usepackage{tikz}
\usepackage[all]{xy}
\usepackage[unicode=true]{hyperref}

\newtheorem{theorem}{Theorem}[section]
\newtheorem{corollary}[theorem]{Corollary}
\newtheorem{criterion}[theorem]{Criterion}
\newtheorem{definition}[theorem]{Definition}
\newtheorem{example}[theorem]{Example}

\newtheorem{lemma}[theorem]{Lemma}
\newtheorem{remark}[theorem]{Remark}
\newtheorem{proposition}[theorem]{Proposition}
\def\Ann{\textup{Ann}}
\def\Im{\textup{Im}}
\def\Ker{\textup{Ker}}
\def\Coker{\textup{Coker}}

\def\Spec{\textup{Spec\,}}
\def\Max{\textup{Max\,}}
\def\Pic{\textup{Pic\,}}

\newcommand{\cref}[1]{eq. (\ref{#1})}

\begin{document}
\title[Matrix invertible extensions]{Matrix invertible
extensions over commutative rings. Part III: Hermite rings}
\author{Grigore C\u{a}lug\u{a}reanu, Horia F. Pop, Adrian Vasiu}

\begin{abstract}
We reobtain and often refine prior criteria due to Kaplansky, McGovern, Roitman, Shchedryk, Wiegand, and Zabavsky--Bilavska and obtain new criteria for a Hermite ring to be an \textsl{EDR}. We mention three criteria: (1) a Hermite ring $R$ is an \textsl{EDR} iff for all pairs $(a,c)\in R^2$, the product homomorphism $U(R/Rac)\times U\bigl(R/Rc(1-a)\bigr)\to U(R/Rc)$ between groups of units is surjective; (2) a reduced Hermite ring $R$ is an \textsl{EDR} iff it is a pre-Schreier ring and for each $a\in R$, every zero determinant unimodular $2\times 2$ matrix with entries in $R/Ra$ lifts to a zero determinant matrix with entries in $R$; (3) a B\'{e}zout domain $R$ is an \textsl{EDD} iff for all triples $(a,b,c)\in R^3$ there exists a unimodular pair $(e,f)\in R^2$ such that $(a,e)$ and $(be+af,1-a-bc)$ are unimodular pairs. We use these criteria to show that each B\'{e}zout ring $R$ that is an $(SU)_2$ ring (as introduced by Lorenzini) such that for each nonzero $a\in R$ there exists no nontrivial self-dual projective $R/Ra$-module of rank $1$ generated by $2$ elements (e.g., all its elements are squares), is an \textsl{EDR}.
\end{abstract}

\maketitle

\section{Introduction}

Let $R$ be a commutative ring with identity element $1$. Let $U(R)$, $N(R)$, $J(R)$, $Z(R)$ and $Pic(R)$ be its group of units, its nilpotent radical, its Jacobson radical, its set of zero divisors, and its Picard group (respectively). For $a\in R$, let $\Ann_R(a)$ be its annihilator in $R$. For $n\in\mathbb{N}=\{1,2,\ldots\}$, let $\mathbb{M}_n(R)$ be the $R$-algebra of $n\times n$ matrices with entries in $R$. Let ${GL}_{n}(R)$ be the general linear group of units of $\mathbb M_n(R)$. Let ${SL}_{n}(R):=\{M\in GL_n(R)|\det(M)=1\}$ be the special linear subgroup of ${GL}_{n}(R)$. Let $Um(F)$ be the set of unimodular elements of a free $R$-module $F$; so we have an identity $Um(R)=U(R)$ of sets. We say that $A,B\in\mathbb M_2(R)$ are congruent modulo an ideal $I$ of $R$ if $A-B\in\mathbb M_2(I)$.

In this paper we study B\'{e}zout rings, i.e., rings whose finitely generated ideals are principal. Each B\'{e}zout ring is an arithmetical ring, i.e., its lattice of ideals is distributive and all its localizations at prime ideals are valuation rings (see \cite{jen}, Thms.\ 1 and 2). The ring $R$ is a B\'{e}zout ring iff each diagonal matrix with entries in $R$ admits diagonal reduction (see \cite{LLS}, Thm.\ (3.1)). 

Recall that: (i) $R$ is a Hermite ring in the sense of Kaplansky if $R^2=RUm(R^2)$, equivalently
if each $1\times 2$ matrix with entries in $R$ admits diagonal reduction; (ii) $R$ is an elementary divisor ring, to be abbreviated as EDR, if for each $(m,n)\in\mathbb N^2$ with $m\le n$, every matrix of size either $m\times n$ or $n\times m$ with entries in $R$ admits diagonal reduction, i.e., is equivalent to a matrix whose off diagonal entries are $0$ and whose diagonal entries $a_{1,1},\ldots,a_{m,m}$ are such that $a_{i,i}$ divides $a_{i+1,i+1}$ for all $i\in \{1,\ldots ,m-1\}$; (iii) $R$ is a pre-Schreier ring if each $x\in R$ is primal, i.e., if $x$ divides $yz$, with $(y,z)\in R^2$, then there exists $(u,v)\in R^2$ such that $x=uv$, $u$ divides $y$ and $v$ divides $z$. 

Each Hermite ring $R$ is a B\'{e}zout ring but the converse does not hold (see \cite{GH1}, Ex.\ 3.4, \cite{WW}, Ex.\ 3.3, or \cite{car}, Prop.\ 8). However, a Hermite domain is the same as a B\'{e}zout domain. B\'{e}zout domains are GCD (greatest common divisors exist) domains and hence pre-Schreier domains with trivial Picard groups. Clearly, each \textsl{EDR} is a Hermite ring. 

We use the tools of Parts I and II in order to obtain necessary and sufficient criteria for a Hermite ring (or domain) to be an \textsl{EDR} (or an \textsl{EDD}, i.e., an elementary divisor domain).

Let $A\in Um\bigl(\mathbb M_2(R)\bigr)$. Recall that $A$ is called extendable if it is obtained from a matrix $A^+\in SL_3(R)$ by removing its third row and its third column (see \cite{CPV1}, Def.\ 1.1). If we can choose $A^+$ such that its $(3,3)$ entry is $0$, then $A$ is called simply extendable. Recall that $A$ is called determinant liftable (resp.\ weakly determinant liftable) if there exists $B\in Um\bigl(\mathbb M_2(R)\bigr)$ (resp.\ $B\in \mathbb M_2(R)$) congruent to $A$ modulo $R\det(A)$ and with $\det(B)=0$ (see \cite{CPV2}, Def.\ 1.1). The {\it only} general implications between these 4 notions on $A$ are the ones in the diagram
\begin{equation}\label{EQ1}
\xymatrix@R=10pt@C=21pt@L=2pt{
\textup{simply extendable} \ar@{=>}[r]\ar@{=>}[d] & \textup{extendable}\ar@{=>}[d]\\
\textup{determinant liftable} \ar@{=>}[r] & \textup{weakly determinant liftable}\\}
\end{equation}
(see \cite{CPV2}, Thm.\ 1.3 for the vertical ones; for the fact that their converses do not hold and that there exist no implications between extendable and determinant liftable see \cite{CPV1}, Exs.\ 5.3 and 6.1 and \cite{CPV2}, Exs.\ 5.1 and 1.9). Recall that a matrix $B\in\mathbb M_2(R)$ is called non-full if it is the product of $2\times 1$ with $1\times 2$ matrices with entries in $R$. An integral domain $R$ is a pre-Schreier domain iff each matrix in $\mathbb{M}_{2}(R)$ of zero determinant is non-full by \cite{MR}, Lem.\ 1.

Recall from \cite{CPV1}, Def.\ 1.2 that $R$ is called a $\Pi_2$ ring if each matrix in $Um\bigl(\mathbb M_2(R)\bigr)$ of zero determinant is extendable, equivalently it is simply extendable by \cite{CPV1}, Lem.\ 4.1(1), and that $R$ is called an $E_2$ (resp.\ an $SE_2$) ring if each matrix in $Um\bigl(\mathbb M_2(R)\bigr)$ is extendable (resp.\ simply extendable). If $Pic(R)$ is trivial (say, $R$ is a B\'{e}zout domain), then $R$ is a $\Pi_2$ ring (see \cite{CPV1}, paragraph after Thm.\ 1.4). 

See \cite{CPV1}, Def. 1.5 for stable ranges notation. If $R$ is a Hermite ring, then $sr(R)\le 2$ (see \cite{MM}, Prop.\ 8(i) or \cite{zab1}, Cor.\ 2.1.1). In fact, a B\'{e}zout ring $R$ is a Hermite ring iff $sr(R)\le 2$ (see \cite{zab1}, Thm.\ 2.1.2). B\'{e}zout domains that have stable range 1 or 1.5 were also studied in \cite{rus} and respectively \cite{shc3}, \cite{shc4} and \cite{BS}. 

For Hermite rings, Parts I and II get summarized in an example and a theorem as follows. 

\begin{example}\normalfont\label{EX1} {\bf (1)} For Hermite rings $R$, simply extendable and extendable properties on a matrix in $Um\bigl(\mathbb M_2(R)\bigr)$ are equivalent (see \cite{CPV1}, Thm.\ 1.6). 

\smallskip
{\bf (2)} For Hermite rings $R$ that are $\Pi_2$ rings, simply extendable, extendable and determinant liftable properties on a matrix in $Um\bigl(\mathbb M_2(R)\bigr)$ are equivalent (see \cite{CPV2}, Cor.\ 1.5). 

\smallskip
{\bf (3)} For rings $R$ such that each zero determinant matrix in $\mathbb M_2(R)$ is non-full (e.g., if $R$ is a product of pre-Schreier domains), extendable and weakly determinant liftable properties on a matrix in $Um\bigl(\mathbb M_2(R)\bigr)$ are equivalent (see \cite{CPV2}, Thm.\ 1.7). 

\smallskip
{\bf (4)} From parts (1) to (3) it follows that for Hermite rings $R$ with the property that zero determinant matrices in $\mathbb M_2(R)$ are non-full, the simply extendable, extendable, determinant liftable, and weakly determinant liftable properties on a matrix in $Um\bigl(\mathbb M_2(R)\bigr)$ are equivalent (cf.\ \cite{CPV2}, Cor.\ 1.8).
\end{example}

\begin{theorem}\label{TH1} Let $R$ be a Hermite ring. Then the following statements are equivalent. 

\medskip \textbf{(1)} The ring $R$ is an \textsl{EDR}.

\smallskip \textbf{(2)} The ring $R$ is an $SE_2$ ring.

\smallskip \textbf{(3)} The ring $R$ is an $E_2$ ring.

\smallskip \textbf{(4)} For each $a\in R$, $R/Ra$ is a $\Pi_2$ ring (equivalently, every projective $R/Ra$-module of rank $1$ generated by $2$ elements is free).

\smallskip \textbf{(5)} Each matrix in $Um\bigl(\mathbb M_2(R)\bigr)$ is determinant liftable and $R$ is a $\Pi_2$ ring.

\medskip
If all zero determinant matrices in $\mathbb M_2(R)$ are non-full, then these five statements are also equivalent to:

\medskip \textbf{(6)} Each matrix in $Um\bigl(\mathbb M_2(R)\bigr)$ is weakly determinant liftable.
\end{theorem}

Referring to Theorem \ref{TH1}, see \cite{CPV1}, Cor.\ 1.8 for $(1)\Leftrightarrow (2)\Leftrightarrow (3)$, see \cite{CPV1}, Cor.\ 4.2 for $(3)\Leftrightarrow (4)$, see \cite{CPV1}, Thm.\ 1.4\ for the equivalence of statement $(4)$, and see \cite{CPV2}, Thm.\ 1.4 for $(2)\Leftrightarrow (5)$. See \cite{CPV2}, Thm.\ 1.7 for $(3)\Leftrightarrow (6)$. 

The equivalence $(1)\Leftrightarrow (4)$ of Theorem \ref{TH1} can be viewed as a more practical way to check that a Hermite ring is an \textsl{EDR} than \cite{WW}, Thm.\ 2.1 (which does not restrict to only 2 generators) or \cite{ZB}, Thms.\ 3 and 5 (which work with finitely generated projective modules over all quotients of $R$).

Recall that $R$ is a Hermite ring iff for all $(m,n)\in \mathbb{N}^2$ and each $m\times n$ matrix $B$ with entries in $R$, there exists $(M,N)\in {GL}_m(R)\times {GL}_n(R)$ such that $MB$ and $BN$ are both lower (equivalently, upper) triangular (see \cite{GH2}, Thm. 3). Thus, in Theorem \ref{TH1}(4) or (5) it suffices to consider only upper triangular matrices, and it is part of the goals of this paper to show that in fact one can restrict to specific upper triangular matrices (e.g., see Propositions \ref{P1} and \ref{P2}).

\begin{example}\normalfont\label{EX2} For $(a,b,c)\in Um(R^3)$, let $A:=\left[ 
\begin{array}{cc}
a & b \\ 
0 & c\end{array}\right]\in Um\bigl(\mathbb M_2(R)\bigr)$.

\medskip
{\bf (1)} The matrix $A$ is determinant liftable iff there exists $(x,y,z,w)\in R^4$ such that $ax+by+cw=1$ and $xw=yz$ (see \cite{CPV2}, Thm.\ 1.2). 

\smallskip
{\bf (2)} If $R$ is arbitrary (resp.\ if $N(R)=0$), then the matrix $A$ is weakly determinant if (resp.\ iff) there exists $(x,y,z,w)\in R^4$ such that $1-ax-by-cw+ac(xw-yz)=0$ by Thm.\ 6.1(2) (resp.\ by Thm.\ 6.1(2) and (3)). As $1-ax-by-cw+ac(xw-yz)$ can be rewritten as $(1-ax)(1-cw)-y(b+acz)$ (cf.\ \cite{CPV2}, Eq.\ (I)), it follows that $A$ is weakly determinant liftable if (resp.\ iff) there exists $(x,z,w)\in R^3$ such that $b+acz$ divides $(1-ax)(1-cw)$.
\end{example}

Example \ref{EX2}(2) and the equivalence $(1)\Leftrightarrow (5)$ of Theorem \ref{TH1} (with (5) stated for upper triangular matrices in $Um\bigl(\mathbb M_2(R)\bigr)$) imply directly the following result.

\begin{corollary}\label{C1}
Assume $R$ is a Hermite ring and a $\Pi_2$ ring. Then $R$ is an \textsl{EDR} iff for each $(a,b,c)\in Um(R^3)$ there exists $(x,y,z,w)\in R^4$ such that $ax+by+cw=1$ and $xw=yz$ (equivalently, there exists $(x,y,w)\in R^3$ such that $ax+by+cw=1$ and $xw\in Ry$).
\end{corollary}

For Hermite rings, to get necessary and/or sufficient conditions to be \textsl{EDR}s that involve units, one would like to get unit interpretations that would characterize the inequality $sr(R)\le 2$ in some way similar to \cite{CPV1}, Prop.\ 2.4. As a progress in this direction we introduce using units the following class of rings.

\begin{definition}\label{def1}
We say that $R$ is a $U_2$ ring if for each $(a,b)\in Um(R^{2})$ and $c\in R$, the natural product
homomorphism 
\begin{equation*}
\upsilon_{a,b,c}:U(R/Rac)\times U(R/Rbc)\rightarrow U(R/Rc)
\end{equation*}is surjective, i.e., the functorial commutative diagram
\begin{equation}\label{EQ2}
\xymatrix@R=10pt@C=21pt@L=2pt{
U\bigl(R/(Rac\cap Rbc)\bigr) \ar[r]^{\;\;\;} \ar[d] & U(R/Rac)\ar[d]\\
U(R/Rbc) \ar[r]^{\;\;\;} & U(R/Rc),\\}
\end{equation}
whose arrows are natural reductions, is not only a pullback but also a pushout.
\end{definition}

For $(a,b)\in Um(R^{2})$, if $(a^{\prime},b^{\prime})\in R^2$ is such that $aa^{\prime}+bb^{\prime}=1$, then $\upsilon_{(aa^{\prime},1-aa^{\prime},c)}$ factors through $\upsilon_{a,b,c}$. Hence $R$ is a $U_2$ ring iff $\upsilon_{a,1-a,c}$ is surjective for all $(a,c)\in R^2$. The fact that a unit in $U(R/Rc)$ does or does not belong to $\Im(\upsilon_{a,1-a,c})$ relates to a certain zero determinant matrix being or not being non-full (see Example \ref{EX10}).

\begin{example}\normalfont\label{EX3}
Assume $asr(R)=1$, i.e., $R$ has almost stable range $1$. We show that $R$ is a $U_2$ ring. Let $(a,b)\in Um(R^{2})$ and $c\in R$. If $c\in J(R)$, then $U(R)$ surjects onto $U(R/Rc)$. If $c\notin J(R)$, then, as $c$ is a linear combination of $ac$ and $bc$, either $ac\notin J(R)$ or $bc\notin J(R)$. Hence either $U(R/Rac)\rightarrow U\bigl(R/(Rac+Rc)\bigr)=U(R/Rc)$ or
$U(R/Rbc)\rightarrow U\bigl(R/(Rbc+Rc)\bigr)=U(R/Rc)$ is a surjective homomorphisms 
(see \cite{CPV1}, Prop.\ 2.4(3)). We incur that $\upsilon_{a,b,c}$ is surjective. So $R$ is a $U_2$ ring.
\end{example}

The next two examples underline the relevance of pre-Schreier rings.

\begin{example}
\normalfont\label{EX4} Let $R$ be such that each zero determinant matrix $B\in\mathbb M_2(R)$ is non-full. We show that $R$ is a pre-Schreier ring. Let $(x,y,z)\in R^3$ be such that $x$ divides $yz$. Let $w\in R$ be such that $xw=yz$. Let $C:=\left[ 
\begin{array}{cc}
x & y \\ 
z & w\end{array}\right]\in \mathbb M_2(R)$. As $\det(C)=0$, there exists $(l,m,o,q)\in R^4$ such that $C=\left[ 
\begin{array}{c}
l \\ 
m\end{array}\right] \left[ 
\begin{array}{cc}
o & q\end{array}\right] $. Hence $x=lo$ with $l$ dividing $y=lq$ and $o$ dividing $z=mo$. Thus $R$ is a pre-Schreier ring.
\end{example}

\begin{example}
\normalfont\label{EX5} Assume $N(R)=0$. Let $B\in\mathbb M_2(R)$ with $\det(B)=0$. We check that if $B$ admits diagonal reduction, then $B$ is non-full. Let $M,N\in GL_2(R)$ be such that $MBN=\left[ 
\begin{array}{cc}
a & 0 \\ 
0 & ab\end{array}\right]$ with $(a,b)\in R^2$. As 
$$a^2b=\det(MBN)=\det(M)\det(B)\det(N)=0,$$ it follows that $(ab)^2=0$ and thus, as $N(R)=0$, we have $ab=0$. Therefore, $B=M^{-1}\left[ 
\begin{array}{cc}
a & 0 \\ 
0 & 0\end{array}\right]N^{-1}=M^{-1}\left[ 
\begin{array}{c}
a \\ 
0\end{array}\right] \left[ 
\begin{array}{cc}
1 & 0\end{array}\right]N^{-1}$ is non-full.
\end{example}

Examples \ref{EX4} and \ref{EX5} imply directly the following result.

\begin{corollary}\label{C1.5}
If $R$ is an \textsl{EDR} with $N(R)=0$, then $R$ is a pre-Schreier ring.
\end{corollary}

Section \ref{S3} proves the following theorem. 

\begin{theorem} 
\label{TH2} We consider the following conditions.

\medskip \textbf{(1)} The ring $R$ is a pre-Schreier ring and for each $(a,b,c)\in Um(R^3)$ there exists $(x,y,z,w)\in R^4$ such that $1-ax-by-cw+ac(xw-yz)=0$.

\smallskip \textbf{(2)} The ring $R$ is a $\Pi_2$ ring and each upper triangular matrix in $Um\bigl(\mathbb M_2(R)\bigr)$ is simply extendable (e.g., $R$ is an $SE_2$ ring).

\smallskip \textbf{(3)} Each upper triangular matrix in $Um\bigl(\mathbb M_2(R)\bigr)$ is extendable and $R$ is an integral domain.

\medskip
 If any one of these conditions holds, then $R$ is a $U_2$ ring.
\end{theorem}

\begin{example}\normalfont\label{EX6}
Assume each upper triangular matrix in $Um\bigl(\mathbb M_2(R)\bigr)$ is weakly determinant liftable (e.g., this holds if $R$ is a $WJ_{2,1}$ ring in the sense of \cite{CPV2}, Def.\ 1.10(1), see \cite{CPV2}, Thm.\ 1.11(1)). If $R$ is also a pre-Schreier ring with $N(R)=0$, then condition (1) of Theorem \ref{TH2} holds by Example \ref{EX2}(2), hence $R$ is a $U_2$ ring.
\end{example}

Each Dedekind domain is a $U_2$ ring by Example \ref{EX3}. Hence Dedekind domains which are not principal ideal domains are $U_2$ rings which (by \cite{CPV1}, Thm.\ 1.7(4)) are not $\Pi_2$ rings. 

In Section \ref{S4} we prove the following `units supplement' to Theorem \ref{TH1}.

\begin{theorem}
\label{TH3} For a Hermite ring $R$ the following statements are equivalent.

\medskip \textbf{(1)} The ring $R$ is an \textsl{EDR}.

\smallskip \textbf{(2)} The ring $R/N(R)$ is a pre-Schreier ring and each matrix in $Um\bigl(\mathbb M_2(R)\bigr)$ is weakly determinant liftable.

\smallskip \textbf{(3)} The ring $R$ is a $U_{2}$ ring.

\smallskip \textbf{(4)} Given $\bigl((a,b),(c,d)\bigr)\in Um(R^2)^2$, there exists $t\in R$ such that we can factor $d+ct=d_1d_2$ with $\bigl((a,d_1),(b,d_2)\bigr)\in Um(R^2)^2$.

\smallskip \textbf{(5)} Given $(a,d)\in R^2$ and $c\in 1+Rd$, 
there exists $t\in R$ such that we can factor $d+ct=d_1d_2$ with $\bigl((a,d_1),(1-a,d_2)\bigr)\in Um(R^2)^2$.\end{theorem}

\begin{example}
\normalfont\label{EX7} Let $R$ be a Hermite ring which is not an \textsl{EDR} (see \cite{GH1}, Sect.\ 4 and Ex.\ 4.11 and \cite{car}, Prop.\ 8). From Theorems \ref{TH1} and \ref{TH3} it follows that $R$ is neither an $E_{2}$ ring nor a $U_{2}$ ring and there exists $a\in R$ such that the Hermite ring $R/Ra$ is not a $\Pi_2$ ring.
\end{example} 

Theorem \ref{TH3} implies directly the following result.\footnote{Corollaries \ref{C1} and \ref{C2} can be also deduced from Example \ref{EX2}(2) and Theorem \ref{TH1} via the following result proved in \cite{GV}. Let $S$ be a Bézout ring $S$ with $N(S)=0$. Then $S$ is a pre-Schreier ring iff all zero determinant matrices in $\mathbb M_2(S)$ are non-full.}

\begin{corollary}\label{C2}
Assume $R$ is a Hermite ring and a pre-Schreier ring such that $N(R)=0$. Then $R$ is an \textsl{EDR} iff each matrix in $Um\bigl(\mathbb M_2(R)\bigr)$ is weakly determinant liftable and iff for each $(a,b,c)\in Um(R^3)$ there exists $(x,y,z,w)\in R^4$ such that
\begin{equation}\label{EQ3}
(1-ax)(1-cw)=y(b+acz).
\end{equation}
\end{corollary}

If $asr(R)=1$, then $R$ is an $U_2$ ring (see Example \ref{EX3}); so McGovern's theorem that each Hermite ring $R$ with $asr(R)=1$ is an \textsl{EDR}  (see \cite{mcg1}, Thm.\ 3.7 or \cite{CPV1}, Cor.\ 1.9(2)) also follows from the equivalence $(1)\Leftrightarrow (3)$ of Theorem \ref{TH3}.

In Section \ref{S5} we prove the following `$\Pi_2$ supplement' to Theorem \ref{TH1}.

\begin{theorem}
\label{TH4} Assume $R$ is a Hermite ring with $N(R)=0$. Then $R$ is an \textsl{EDR} iff for each $a\in R$ the quotient ring $R/\Ann_R(a)$ is a $\Pi_2$ ring and each matrix in $Um\bigl(\mathbb M_2(R)\bigr)$ is weakly determinant liftable (equivalently, and for every $(a,b,c)\in Um(R^3)$ there exists $(x,y,z,w)\in R^4$ such that Equation (\ref{EQ3}) holds).
\end{theorem}

Lorenzini introduced 3 classes of rings that are `between' Hermite rings and \textsl{EDR}s (see \cite{lor}, Prop.\ 4.11). For the first class, called $J_{2,1}$ (see \cite{lor}, Def.\ 4.6), it was proved in \cite{CPV2}, Thm.\ 1.11(2) that each $J_{2,1}$ ring which is a $\Pi_2$ ring is an \textsl{EDR}. We recall the last 2 classes of rings (see \cite{lor}, Def.\ 2.1) in a slightly different way.

\begin{definition}\label{def2}
For $n\ge 2$ we say that $R$ is:

\medskip \textbf{(1)} a ${(WSU^{\prime})}_{n}$ (resp.\ ${(WSU)}_{n}$) ring if for each $A\in Um\bigl(\mathbb{M}_n(R)\bigr)$ there exists $N\in {SL}_n(R)$ (resp.\ $N\in {GL}_n(R)$) such that $AN$ is symmetric;

\smallskip \textbf{(2)} an ${(SU^{\prime})}_{n}$ (resp.\ ${(SU)}_{n}$) ring if it is a Hermite ring and a ${(WSU^{\prime})}_{n}$ (resp.\ ${(WSU)}_{n}$) ring.
\end{definition}

In Definition \ref{def2}(1), the existence of $N\in {SL}_n(R)$ (resp.\ $N\in {GL}_n(R)$) is equivalent to the existence of $M\in {SL}_n(R)$ (resp.\ $M\in {GL}_n(R)$) such that $MA$ is
symmetric and is equivalent to the existence of a pair $(M,N)\in {SL}_n(R)^2$ (resp.\ $(M,N)\in {GL}_n(R)^2$) such that $MAN$ is symmetric. This follows via conjugation from the fact that for each symmetric matrix $O\in\mathbb M_n(R)$ and every $M\in {GL}_n(R)$, $MOM^T$ is symmetric.

Each $(WSU^{\prime})_n$ ring is a $(WSU)_n$ ring. As the extendable (resp.\ simply extendable) property of a unimodular $2\times 2$ matrix depends only on the equivalence class of the matrix (see \cite{CPV1}, Lem.\ 4.1(3)), it follows that a $(WSU)_2$ ring is an $E_2$ (resp.\ $SE_2$ ring) iff each symmetric matrix in $Um\bigl(\mathbb M_2(R)\bigr)$ is extendable (resp.\ simply extendable). Section \ref{S6} proves basic properties of ${(WSU^{\prime})}_n$ and ${(WSU)}_n$ rings; e.g., our definition of ${(SU^{\prime})}_{n}$ or ${(SU)}_{n}$ rings is equivalent to the one in \cite{lor} by Proposition \ref{P4}. 

In Section \ref{S7} we prove the following theorem.

\begin{theorem}
\label{TH5} Let $R$ be a $(WSU)_2$ ring. Then the following properties hold.

\medskip \textbf{(1)} For each $(a,b)\in Um(R^2)$ and $c\in R$, $\Coker(\upsilon_{a,b,c})$ is a Boolean group, i.e., we have an inclusion $\{x^2|x\in U(R/Rc)\}\subset\Im(\upsilon_{a,b,c})$.

\smallskip \textbf{(2)} If $sr(R)\le 4$, then for each $d\in R$, every projective $R/Rd$-module of rank $1$ generated by $2$ elements is self-dual.

\smallskip \textbf{(3)} Assume that each element of $R$ is the square of an element of $R$ (e.g., this holds if $R$ is an integrally closed domain with an algebraically closed field of fractions or is a perfect ring of characteristic $2$). Then $R$ is a $U_{2}$ ring. If moreover $R$ is a Hermite ring, then $R$ is an \textsl{EDR}.
\end{theorem}

In Section \ref{S8} we prove the following theorem.

\begin{theorem}
\label{TH6}
Let $R$ be a Hermite ring such that for each $(a,b)\in Um(R^2)$ and $c\in R$, $\Coker(\upsilon_{a,b,c})$ is a Boolean group. Then for each $d\in R$, every projective $R/Rd$-module of rank $1$ generated by $2$ elements is self-dual.
\end{theorem}

To connect with Pell-type equations and provide more examples of $(WSU)_2$ rings that are \textsl{EDR}s, we first prove in Section \ref{S9} the following non-full Pell-type criterion.

\begin{criterion}
\label{CR1} Let $A=\left[ 
\begin{array}{cc}
a & b \\ 
b & c\end{array}\right]\in Um\bigl(\mathbb{M}_2(R)\bigr)$ be symmetric and with zero determinant. Then the
following properties hold.

\medskip \textbf{(1)} The matrix $A$ is simply extendable if there exists $(e,f)\in R^{2}$ such that $ae^{2}-cf^{2}\in U(R)$, and the converse holds if 
$R$ has characteristic $2$.

\smallskip \textbf{(2)} Assume that $R$ is a Hermite ring of characteristic $2$ with $N(R)=0$. If $b$ is not a zero divisor, then $A$ is simply extendable.
\end{criterion}

By combining Theorem \ref{TH5}(1) with Criterion \ref{CR1}, we obtain the
following Pell-type criterion proved in Section \ref{S9}.

\begin{criterion}
\label{CR2} Let $R$ be an $(SU)_2$ ring. Then $R$ is an \textsl{EDR} if for all triples $(a,b,c)\in Um(R^3)$ there exists $(e,f)\in R^2$ such that $(ae^2-cf^2,ac-b^2)\in Um(R^2)$, and the converse holds if $R$ has
characteristic $2$. 
\end{criterion}

We recall that, if the characteristic of $R$ if a prime, then its perfection $R_{{perf}}$ is the inductive limit of the inductive system indexed by $n\in\mathbb{N}$ whose all transition homomorphisms are the Frobenius endomorphism of $R$. If $R$ is a Hermite (or \ B\'{e}zout) ring, then so is $R_{{perf}}\cong \bigl(R/N(R)\bigr)_{{perf}}$.

\begin{example}
\normalfont\label{EX9} Assume $R$ has characteristic $2$. If for all $(a,c)\in R^2$, $\Coker(\upsilon_{a,1-a,c})$ is a Boolean group, then $R_{{perf}}$ is a $U_{2}$ ring. From this and Theorem \ref{TH5}(3) it follows that perfections of $(WSU)_{2}$ rings of characteristic $2$ are $U_{2}$ rings. Thus perfections of $(SU)_{2}$ rings of characteristic $2$ are \textsl{EDR}s by (the implication $(3)\Rightarrow (1)$ of) Theorem \ref{TH3}.
\end{example}

As an application of Corollary \ref{C2} in Section \ref{S10} we prove the following criterion.

\begin{criterion}
\label{CR3} For a B\'{e}zout domain $R$ the following statements are equivalent.

\medskip \textbf{(1)} The ring $R$ is an \textsl{EDD}.

\smallskip \textbf{(2)} For each $(a,b,s)\in R^3$, there exists a pair $(q,r)\in R^2$ such that by defining $y:=r+s-asq-bqr$ and $t:=1+q-aq-br$
we have $t\in Ry+Rat$ (equivalently, there exists $(y_1,y_2)\in R^2$ such that $y=y_1y_2$, $y_1$ divides $t$, and $(a,y_2)\in Um(R^2)$).

\smallskip \textbf{(3)} For each $(a,b,s)\in R^3$ there exists a pair $(e,f)\in Um(R^2)$ such that $\bigl((a,e),(be+af,1-bs-a)\bigr)\in Um(R^2)^2$.
\end{criterion}

E.g., statement (3) of Criterion \ref{CR3} holds if $(a,s)\in Um(R^2)$ as we can take $(e,f):=(s,1)$, if $(1-a,b)\in Um(R^2)$ as we can take $(e,f):=(1,0)$, or if there exists $q\in R$ such that $(b+aq,1-bs-a)\in Um(R^2)$ as we can take $(e,f):=(1-a,q+b)$. Hence directly from Criterion \ref{CR3} we get the following consequence.

\begin{corollary}
Assume $R$ is a B\'{e}zout domain with the property that for all $(a,b,s)\in R^3$ with $\bigl((a,s),(b,1-a)\bigr)\notin Um(R^2)^2$ there exists $q\in R$ such that $(b+aq,1-bs-a)\in Um(R^2)$. Then $R$ is an \textsl{EDD}.
\end{corollary}

The implicit and explicit questions raised in the literature, such as,
\textquotedblleft Is a B\'{e}zout domain of finite Krull dimension [at least 
$2$] an \textsl{EDD}?" (see \cite{FS}, Ch.\ III, Probl.\ 5, p.\ 122), and,
`What classes of B\'{e}zout domains which are not \textsl{EDD}s exist?',
remain unanswered. However, the above results reobtain or can be easily used to
reobtain multiple other criteria existing in the literature of when a Hermite 
ring is an \textsl{EDR}. Moreover, the equivalence $(1)\Leftrightarrow (4)$ of Theorem \ref{TH3} was proved for B\'{e}zout domains in \cite{shc2}, Thm.\ 5 and a variant of it was proved for Hermite rings in \cite{rot}, Prop.\ 2.9. Equation (\ref{EQ3}) generalizes and refines the equation one would get based on \cite{rot}, Rm.\ 2.8. The last two references were reinterpreted in terms of neat stable range $1$ (see \cite{zab2}, Def.\ 21) in \cite{zab2}, Thms.\ 31 and 33; see \cite{mcg2} for clean and neat rings.

\section{Test matrices}\label{S2}

In this section we introduce companion and universal test matrices for unimodular $2\times 2$ matrices and prove a few properties of them that will be often used in the subsequent sections to streamline proofs. We begin with companion test matrices.

\begin{definition}\label{def3}
Let $A\in Um\bigl(\mathbb M^2(R)\bigr)$. An upper triangular matrix $B\in Um\bigl(\mathbb M^2(R)\bigr)$ is called a companion test matrix for $A$ if there exists $(a,b,c,a^{\prime},c^{\prime})\in R^5$ such that $A$ is equivalent to $\left[\begin{array}{cc}
a & b \\ 
0 & c\end{array}\right] $ and $B=\left[\begin{array}{cc}
aa^{\prime} & b \\ 
0 & cc^{\prime}\end{array}\right]$. 
\end{definition}

Let $\mathcal P$ be one of the 4 notions extendable, simply extendable, determinant liftable and weakly determinant liftable. Definition \ref{def3} is justified by the next proposition. 

\begin{proposition}\label{P1}
If $\mathcal P$ is weakly determinant liftable we assume that $N(R)=0$. If $A\in Um\bigl(\mathbb M^2(R)\bigr)$ has a companion test matrix $B$ which is $\mathcal P$, then $A$ is $\mathcal P$.
\end{proposition}

\begin{proof}
As $\mathcal P$ property depends only on equivalence classes (see \cite{CPV1}, Lem.\ 4.1(3) for the case of extendable and simply extendable notions), it suffices to show that if $(a,b,c,a^{\prime},c^{\prime})\in R^5$ is such that $B:=\left[\begin{array}{cc}
aa^{\prime} & b \\ 
0 & cc^{\prime}\end{array}\right]$ is $\mathcal P$, then $C:=\left[\begin{array}{cc}
a & b \\ 
0 & c\end{array}\right]$ is $\mathcal P$. 

First we assume that $\mathcal P$ is extendable (resp.\ simply extendable). Let $(e,g)\in R^2$ (resp.\ $(e,g)\in Um(R^2)$) be such
that $(aa^{\prime}e,be+cc^{\prime}g,aa^{\prime}cc^{\prime})\in Um(R^3)$ (resp.\ $(aa^{\prime}e,be+cc^{\prime}g)\in Um(R^2)$) by \cite{CPV1}, Cor.\ 4.7(1) (resp.\ \cite{CPV1}, Thm.\ 4.3) applied to $B$. Let $f:=c^{\prime}g\in R$. As $(aa^{\prime}e,be+cf,aa^{\prime}cc^{\prime})\in Um(R^3)$ (resp. $(aa^{\prime}e,be+cf)\in Um(R^2)$), it follows that $(ae,be+cf,ac)\in Um(R^3)$ (resp.\ $(ae,be+cf)\in Um(R^{2})$). Thus $C$ is $\mathcal P$ by \cite{CPV1}, Cor.\ 4.7(1) (resp.\ \cite{CPV1}, Thm.\ 4.3).

Next we assume that $\mathcal P$ is determinant liftable (resp.\ weakly determinant liftable). From Example \ref{EX2}(1) (resp.\ \ref{EX2}(2)) applied to $B$ it follows that there exists a quadruple $(x,y,z,w)\in R^4$ such that $1-aa^{\prime}x-by-cc^{\prime}w=0$ and $xw=yz$ (resp.\ such that $1-aa^{\prime}x-by-cc^{\prime}w+aa^{\prime}cc^{\prime}(xw-yz)=0$). For the quadruple $(x^{\prime},y^{\prime},z^{\prime},w^{\prime}):=(a^{\prime}x,y,a^{\prime}c^{\prime}z,c^{\prime}w)\in R^4$, we compute that $1-ax^{\prime}-by^{\prime}-cw^{\prime}=0$ and $x^{\prime} w^{\prime}=y^{\prime}z^{\prime}$ (resp.\ $1-ax^{\prime}-by^{\prime}-cw^{\prime}+ac(x^{\prime}w^{\prime}-y^{\prime}z^{\prime})=0$), so $C$ is $\mathcal P$ by Example \ref{EX2}(1) (resp.\ \ref{EX2}(2)).
\end{proof}

For universal test matrices we consider two cases as follows.

{\bf Case 1: Hermite rings.} In this and the next paragraph we assume that $R$ is a Hermite ring, i.e., for every
pair $(p,q)\in R^{2}$ there exists $(r,s,t)\in R^3$ such that $p=rs$, $q=rt$ and $(s,t)\in Um(R^{2})$. If moreover $R$ is an integral domain (i.e., if $R$ is a B\'{e}zout domain), then $r$ is unique up to a multiplication with a unit of $R$ and is called the greatest common divisor of $p$ and $q$ and one writes $r=\gcd (x,y)$. Our convention is $\gcd (0,0)=0$ (so that we can
still write $0=0\cdot 1$ with $\gcd (1,1)=1$).

Let $A\in Um\bigl(\mathbb{M}_{2}(R)\bigr)$ and $M\in {SL}_2(R)$ be such that 
$B:=MA=\left[ 
\begin{array}{cc}
g & u \\ 
0 & h\end{array}\right] $ is upper triangular. We write $g=ac$ and $h=bc$ with $(a,b)\in
Um(R^{2})$ and $c\in R$. As $B\in Um\bigl(\mathbb{M}_{2}(R)\bigr)$ and $Rg+Rh=Rc$, we
have $(c,u)\in Um(R^{2})$. Let $(a^{\prime },b^{\prime },c^{\prime
},u^{\prime })\in R^4$ be such that $aa^{\prime }+bb^{\prime }=cc^{\prime
}+uu^{\prime }=1$. For $d\in R$, $A$ and $B$ are equivalent to 
$$C:=\left[ 
\begin{array}{cc}
1 & db^{\prime } \\ 
0 & 1\end{array}\right] B\left[ 
\begin{array}{cc}
1 & da^{\prime } \\ 
0 & 1\end{array}\right] =\left[ 
\begin{array}{cc}
ac & u+cd(aa^{\prime }+bb^{\prime }) \\ 
0 & bc\end{array}\right] =\left[ 
\begin{array}{cc}
ac & u+cd \\ 
0 & bc\end{array}\right].$$ 
Here the role of $u+cd$ is that of an arbitrary element of $R$
whose reduction modulo $Rc$ is the `fixed' unit $u+Rc\in U(R/Rc)$. Let $D_{a',b',d}:=\left[ 
\begin{array}{cc}
aa^{\prime }cc^{\prime } & u+cd \\ 
0 & bb^{\prime }cc^{\prime }\end{array}\right]$. If $D_{a',b',d}\in Um\bigl(\mathbb M_2(R)\bigr)$, then $D_{a',b',d}$ is a companion test matrix for $A$ (or $B$). Note that $D_{a',b',0}$ is in $Um\bigl(\mathbb M_2(R)\bigr)$ and is the image of the `first universal test matrix for Hermite
rings' 
\begin{equation*}
\mathcal{D}:=\left[ 
\begin{array}{cc}
x(1-yz) & y \\ 
0 & (1-x)(1-yz)\end{array}\right] \in Um\bigl(\mathbb{M}_{2}(\mathbb{Z}[x,y,z])\bigr)
\end{equation*}via the ring homomorphism $\mathbb{Z}[x,y,z]\rightarrow R$ that maps $x $, $y$, and $z$ to $aa^{\prime }$, $u$, and $u^{\prime }$ (respectively); so $1-x$
maps to $bb^{\prime }=1-aa^{\prime }$ and $1-yz$ maps to $1-uu^{\prime
}=cc^{\prime }$. We also note that $\mathcal{D}$ is equivalent to the matrix 
\begin{equation*}
\mathcal{E}:=\left[ 
\begin{array}{cc}
1 & 0 \\ 
z(x-1)(1-yz) & 1\end{array}\right] \mathcal{D}\left[ 
\begin{array}{cc}
1 & 0 \\ 
xz & 1\end{array}\right] =\left[ 
\begin{array}{cc}
x & y \\ 
0 & (1-x)(1-yz)^{2}\end{array}\right]
\end{equation*}which is the image of the `second universal test matrix for Hermite rings'
\begin{equation*}
\mathcal{F}:=\left[ 
\begin{array}{cc}
x & y \\ 
0 & (1-x)(1-yz)\end{array}\right] \in Um\bigl(\mathbb{M}_{2}(\mathbb{Z}[x,y,z])\bigr)
\end{equation*}via the endomorphism of $\mathbb{Z}[x,y,z]$ that fixes $x$ and $y$ and maps $z$ to $2z-yz^{2}$. See Corollary \ref{C4} for the usage of universal test matrix in this paragraph.

{\bf Case 2: all rings.} The `universal test upper triangular matrix for all rings' is 
\begin{equation*}
\mathcal{G}:=\left[ 
\begin{array}{cc}
x & y \\ 
0 & 1-x-yz\end{array}\right] \in Um\bigl(\mathbb{M}_{2}(\mathbb{Z}[x,y,z])\bigr).
\end{equation*}

\begin{proposition}\label{P2} If $\mathcal P$ is weakly determinant liftable we assume that $N(R)=0$. Then the following properties hold.

\medskip \textbf{(1)} Each upper triangular matrix in $Um\bigl(\mathbb M_2(R)\bigr)$ is $\mathcal P$ iff for each
homomorphism $\phi :\mathbb{Z}[x,y,z]\rightarrow R$, the image of $\mathcal{G}\in Um\bigl(\mathbb{M}_{2}(\mathbb{Z}[x,y,z])\bigr)$ in $Um\bigl(\mathbb{M}_{2}(R)\bigr)$ via $\phi $ is $\mathcal P$.

\smallskip \textbf{(2)} Each zero determinant upper triangular matrix in $Um\bigl(\mathbb M_2(R)\bigr)$ is $\mathcal P$ iff for each
homomorphism $\phi :\mathbb{Z}[x,y,z]\rightarrow R$ with $x(1-x-yz)\in\Ker(\phi)$, the image of $\mathcal{G}\in Um\bigl(\mathbb{M}_{2}(\mathbb{Z}[x,y,z])\bigr)$ in $Um\bigl(\mathbb{M}_{2}(R)\bigr)$ via $\phi $ is $\mathcal P$.
\end{proposition}

\begin{proof}
The `only if' parts are clear. For the `if' part of (1), we consider an upper triangular matrix $A=\left[ 
\begin{array}{cc}
a & b \\ 
0 & c\end{array}\right] \in Um\bigl(\mathbb{M}_{2}(R)\bigr)$. Let $(a^{\prime },b^{\prime },c^{\prime
})\in R^3$ be such that $aa^{\prime }+bb^{\prime }+cc^{\prime }=1$. Let $\phi :\mathbb{Z}[x,y,z]\rightarrow R$ be the homomorphism that maps $x$, $y$ and $z$ to $aa^{\prime }$, $b$ and $b^{\prime }$ (respectively). The image $\left[ 
\begin{array}{cc}
aa^{\prime } & b \\ 
0 & cc^{\prime }\end{array}\right]$ of $\mathcal{G}$ via $\phi$ is $\mathcal P$, so $A$ is $\mathcal P$ by Proposition \ref{P1}. The `if' part of (2) is proved similarly as $ac=0$ implies $aa^{\prime }cc^{\prime }=0$ and $x(1-x-yz)\in\Ker(\phi)$.
\end{proof}

\section{Proof of Theorem \ref{TH2}}\label{S3}

Let $(a,b)\in Um(R^{2})$ and $c\in R$. It suffices to prove that $\upsilon_{a,b,c}$ is surjective if any one of conditions (1) to (3) of Theorem \ref{TH2} holds.

For a fixed unit $\bar{u}\in U(R/Rc)$, let $u\in R$ be such that $\bar{u}=u+Rc$ and consider the matrix $A=\left[ 
\begin{array}{cc}
ac & u \\ 
0 & bc\end{array}\right] \in Um\bigl(\mathbb{M}_{2}(R)\bigr)$.

Assume condition (1) of Theorem \ref{TH2} holds. Thus $R$ is a pre-Schreier ring and there exists $(x,y,z,w)\in R^4$ such that $(1-acx)(1-bcw)=y(u+abc^2z)$ by Example \ref{EX2}(2). As $u+abc^2z$ divides $(1-acx)(1-bcw)$ and $R$ is a pre-Schreier ring, we only use the existence of a nonuple $(x,y,z_1,z_2,w,q,r,s,t)\in R^9$, e.g., $\nu:=(x,y,bcz,0,w,1,1,1,1)\in R^9$, such that $(ac,sq),(bc,rt)\in Um(R^2)$ and we can factor $qru+acz_1+bcz_2=u_{a}u_{b}$, with $u_{a}\in R$ dividing $s-acx$ and $u_{b}\in R$ dividing $t-bcw$. Hence the two units $(q+Rac)^{-1}\cdot (u_{a}+Rac)\in U(R/Rac)$ and $(r+Rbc)^{-1}\cdot (u_{b}+Rbc)\in U(R/Rbc)$ are such that the product of their images in $U(R/Rc)$ is $\bar{u}$. Thus $\upsilon_{a,b,c}$ is surjective.

Assume condition (2) of Theorem \ref{TH2} holds. Thus $R$ is a $\Pi_2$ ring and $A$ is simply extendable. Hence $A$ is determinant liftable by Diagram (\ref{EQ1}). Let the quadruple $(x,y,z,w)\in R^4$ be such that $(1-acx)(1-bcw)=y(u+abc^2z)$ and $xw=yz$ (see Example \ref{EX2}(1)). If there exists $\mathfrak m\in\Max R$ which contains $1-acx$, $1-bcw$ and $y$, then it follows that $xw+\mathfrak m=yz+\mathfrak m\in R/\mathfrak m$ is both zero and non-zero, a contradiction. Thus $(1-acx,1-bcw,y)\in Um(R^3)$. We only use the existence of a nonuple $(x,y,z_1,z_2,w,q,r,s,t)\in R^9$, e.g., the $\nu$ above, such that 
\begin{equation*}
D:=\left[ 
\begin{array}{cc}
s-acx & qru+acz_1+bcz_2 \\ 
y & t-bcw\end{array}\right] \in Um\bigl(\mathbb{M}_{2}(R)\bigr)
\end{equation*}has zero determinant, $(ac,sq),(bc,rt)\in Um(R^2)$ and $(s-acx,t-bcw,y)\in Um(R^3)$. As $R$ is a $\Pi _{2}$ ring, $D$ is non-full (see \cite{CPV1}, Thm.\ 1.4) and therefore we can factor $qru+acz_1+bcz_2=u_{a}u_{b}$, with $u_{a}\in R$ dividing $s-acx$ and $u_{b}\in R$ dividing $t-bcw$. As in the previous paragraph we argue that $\upsilon_{a,b,c}$ is surjective.

Assume condition (3) of Theorem \ref{TH2} holds. Thus $R$ is an integral domain and $A$ is extendable. To prove that $\upsilon_{a,b,c}$ is surjective we can assume that $abc\neq 0$. As $A$ is extendable, its
reduction modulo $Rabc^{2}$ is non-full by \cite{CPV1}, Prop.\ 5.1(2). So the system of congruences 
\begin{equation*}
xy\equiv ac\;(\textup{mod}\; abc^{2}),\;\,xw\equiv u\;(\textup{mod}\; abc^{2}),\;\,yz\equiv 0\;(\textup{mod}\; abc^{2}),\;\,zw\equiv bc\;(\textup{mod}\; abc^{2})
\end{equation*}has a solution $(x,y,z,w)\in R^{4}$. As $(c,u)\in Um(R^2)$, from the second congruence it follows
that $(c,xw)\in Um(R^2)$. From this and the first and fourth
congruences, it follows that there exists $(y^{\prime },z^{\prime })\in R^2$ such
that $y=cy^{\prime }$ and $z=cz^{\prime }$. So, as $R$ is an integral domain, the above first system of
congruences is equivalent to a second one
\begin{equation*}
xy^{\prime }\equiv a\;(\textup{mod}\; abc),\;\;xw\equiv u\;(\textup{mod}\; abc^{2}),\;\;y^{\prime
}z^{\prime }\equiv 0\;(\textup{mod}\; ab),\;\;z^{\prime }w\equiv b\;(\textup{mod}\; abc).
\end{equation*}As $(a,b)\in Um(R^{2})$, the last two congruences imply firstly that 
$(a,z^{\prime })\in Um(R^{2})$ and secondly that there exists $y^{\prime
\prime }\in R$ such that $y^{\prime }=ay^{\prime \prime }$. Thus the second
system of congruences is equivalent to a third one 
\begin{equation*}
xy^{\prime \prime }\equiv 1\;(\textup{mod}\; bc),\;\;xw\equiv u\;(\textup{mod}\; abc^{2}),\;\;y^{\prime
\prime }z^{\prime }\equiv 0\;(\textup{mod}\; b),\;\;z^{\prime }w\equiv b\;(\textup{mod}\; abc).
\end{equation*}From the first and third congruences it follows firstly that $(b,y^{\prime
\prime })\in Um(R^{2})$ and secondly that there exists $z^{\prime \prime
}\in R$ such that $z^{\prime }=bz^{\prime \prime }$. Thus the third system
of congruences is equivalent to a fourth one that has only three congruences
\begin{equation*}
xy^{\prime \prime }\equiv 1\;(\textup{mod}\; bc),\;\;xw\equiv u\;(\textup{mod}\; abc^{2}),\;\;z^{\prime
\prime }w\equiv 1\;(\textup{mod}\; ac).
\end{equation*} Thus $x\in
U(R/Rbc)$ and $w+Rac\in U(R/Rac)$ are such that the product of their images in $U(R/Rc)$ is $\bar{u}$. Thus $\upsilon_{a,b,c}$ is surjective. 

Hence Theorem \ref{TH2} holds.

\begin{example}\normalfont\label{EX10} 
For $(a,c,u)\in R^3$, let $A:=\left[ 
\begin{array}{cc}
ac & u \\ 
0 & (1-a)c\end{array}\right]\in\mathbb{M}_{2}(R)$. We assume that $Rac\cap R(1-a)c=0$ and $(c,u)\in Um(R^2)$; so $a(1-a)c=0$ and $A$ is unimodular with zero determinant. We consider the ideals $J_1:=Rac$ and $J_2:=R(1-a)c$ of $R$. For $i\in\{1,2\}$, let $R_i:=R/J_i$, $c_i:=c+J_i\in R_i$, and $u_i:=u+J_i\in R_i$. The reduction of $A$ modulo $J_1$ is $A_1:=\left[\begin{array}{cc}
0 & u_1 \\ 
0 & c_1\end{array}\right]\in Um\bigl(\mathbb{M}_{2}(R_1)\bigr)$ and modulo $J_2$ is $A_2:=\left[\begin{array}{cc}
c_2 & u_2 \\ 
0 & 0\end{array}\right]\in Um\bigl(\mathbb{M}_{2}(R_2)\bigr)$. Both $A_1$ and $A_2$ are non-full. For $i\in\{1,2\}$, to solve the matrix equation $A_i=\left[\begin{array}{c}
l_i\\ 
m_i\end{array}\right]\left[\begin{array}{cc}
o_i & q_i\end{array}\right]$ in $\mathbb M_2(R_i)$ is the same as solving the system of equations
$l_io_i=c_i\delta_{i2},\; l_iq_i=u_i,\; m_io_i=0,\; m_iq_i=c_i\delta_{i1}$ with $(l_i,m_i,o_i,q_i)\in R_i^4$, where $\delta_{i1}$ and $\delta_{i2}$ are Kronecker deltas. As $(c,u)\in Um(R^2)$, we have $q_1\in U(R_1)$ and $l_2\in U(R_2)$. So the solution sets are 
$$\{(u_1q_1^{-1}, c_1q_1^{-1},0,q_1)|q_1\in U(R_1)\}\subset R_1^4$$ 
for $i=1$ and 
$$\{(l_2, 0,c_2l_2^{-1},u_2l_2^{-1})|l_2\in U(R_2)\}\subset R_2^4$$ 
for $i=2$.
As we have an identity $Rc=J_1+J_2$, the matrix $A$ is non-full iff there exists a pair $(q_1,l_2)\in U(R_1)\times U(R_2)$ such that the images of $u_1q_1^{-1}$ and $l_2$ in $R/Rc$ are equal, i.e., $u+Rc=\upsilon_{a,1-a,c}(q_1,l_2)$, and hence iff $u+Rc\in\Im(\upsilon_{a,1-a,c})$.
\end{example}

\section{Proof of Theorem \protect\ref{TH3} and applications}\label{S4}





To prove Theorem \ref{TH3}, let $S:=R/N(R)$. 

We show that $(1)\Rightarrow (2)$. As $R$ is an \textsl{EDR}, each matrix in $Um\bigl(\mathbb M_2(R)\bigr)$ is determinant liftable by (the implication $(1)\Rightarrow (5)$ of) Theorem \ref{TH1} and so it is weakly determinant liftable. Also, from \cite{her}, Thm.\ 3 it follows that $S$ is an \textsl{EDR} and hence it is a pre-Schreier ring by Example \ref{EX5}. So $(1)\Rightarrow (2)$ holds. 

We show that $(2)\Rightarrow (3)$ if $N(R)=0$. As $N(R)=0$, from the `iff' part of Example \ref{EX2}(2) it follows that condition (1) of Theorem \ref{TH2} holds, hence $R$ is a $U_2$ ring. 

We have $(3)\Leftrightarrow (4)$, as for $\bigl((a,b),(c,d)\bigr)\in Um(R^2)^2$, $(4)$ only translates what means that the unit $d+Rc\in U(R/Rc)$ is in the image of $\upsilon_{a,b,c}$. 

Clearly, $(4)\Rightarrow (5)$ by taking $b:=1-a$. 

To check that $(5)\Rightarrow (4)$, let $\bigl((a,b),(c,d)\bigr)\in Um(R^2)^2$. If the quadruple $(e,f,s,t)\in R^4$ is such that $ae+bf=cs+dt=1$, then by applying (5) to $(ae,d)$ and $sc\in 1+Rd$ it follows that there exists $t_1\in R$ such that we can factor $d+t_1sc=d_1d_2$ with $(ae,d_1)\in Um(R_2)$ and $(1-ae,d_2)=(bf,d_2)\in Um(R_2)$ and thus $\bigl((a,d_1),(b,d_2)\bigr)\in Um(R_2)$ and statement (4) holds by taking $t:=t_1s$.

We show that $(3)\Rightarrow (1)$. As $R$ is an \textsl{EDR} iff it is an $E_2$ ring by (the implication $(1)\Rightarrow (3)$ of) Theorem \ref{TH1}, it suffices to show that if $R$ is a $U_{2}$ ring, then each $A\in Um\bigl(\mathbb{M}_{2}(R)\bigr)$ is
extendable. Up to matrix equivalence (see \cite{CPV1}, Lem.\ 4.1(3)), we can assume that $A$ is upper
triangular, and hence we can write $A=\left[ 
\begin{array}{cc}
ac & u \\ 
0 & bc\end{array}\right] $ where $(a,b,c,u)\in R^4$ is such that $\det (A)=abc^{2}$, $(a,b)\in
Um(R^{2})$ and $u$ in $R$ is (up to matrix equivalence, see Case 1 of Section \ref{S2}) an arbitrary representative of a fixed unit $u+Rc\in U(R/Rc)$. As $R$ is a $U_{2}$ ring, based on the arbitrariness part we can assume there exists $(d,e)\in R^2$ such
that $u=de$, $d+Rac\in U(R/Rac)$ and $e+Rbc\in U(R/Rbc)$. Let $d^{\prime },e^{\prime
}\in R$ be such that $d^{\prime }+Rac$ is the inverse of $d+Rac$ and $e^{\prime }+Rac$ is the inverse of $e+Rbc$. Let $(f,g)\in R^2$ be such
that $dd^{\prime }=1+acf$ and $ee^{\prime }=1+bcg$. The non-full matrix 
\begin{equation*}
B:=\left[ 
\begin{array}{cc}
ac(1+bcg) & u \\ 
abc^{2}d^{\prime }e^{\prime } & bc(1+acf)\end{array}\right] =\left[ 
\begin{array}{cc}
acee^{\prime } & de \\ 
abc^{2}d^{\prime }e^{\prime } & bcdd^{\prime }\end{array}\right] \in \mathbb{M}_{2}(R)
\end{equation*}is congruent to $A$ modulo $R\det (A)$. So $A$ modulo $R\det (A)$ is non-full, thus $A$ is extendable by \cite{CPV1}, Prop.\ 5.1(2).

We are left to show that $(2)\Rightarrow (3)$ in general (i.e., without assuming that $N(R)=0$). For $\bar A\in Um\bigl(\mathbb M_2(S)\bigr)$, let $A\in Um\bigl(\mathbb M_2(R)\bigr)$ be such that its reduction modulo $N(R)$ is $\bar A$. As $A$ is weakly determinant liftable, there exists $C\in\mathbb M_2(R)$ congruent to $A$ modulo $R\det(A)$ and $\det(C)=0$. As the reduction $\bar C\in\mathbb M_2(S)$ of $C$ modulo $N(R)$ is congruent to $\bar A$ modulo $S\det(\bar A)$ and has zero determinant, $\bar A$ is weakly determinant liftable. So $S$ has the same properties as $R$ (note that $N(S)=0$ and $S=S/N(S)=R/N(R)$ is a pre-Schreier ring). As we proved that $(2)\Rightarrow (3)$ if $N(R)=0$, it follows that $S$ is a $U_2$ ring. As $(3)\Rightarrow (1)$, we incur that $S$ is an \textsl{EDR}. It is well-known that this implies that $R$ is an \textsl{EDR}: e.g., see \cite{her}, Thm.\ 3; this also follows from Theorem \ref{TH1} via the fact that $A$ is (simply) extendable iff $\bar A$ is so, as one can easily check based on \cite{CPV1}, Thm.\ 4.3 or Cor.\ 4.7. As $R$ is an \textsl{EDR}, it is an $SE_2$ ring by (the implication $(1)\Rightarrow (2)$ of) Theorem \ref{TH1} and hence a $U_2$ ring by Theorem \ref{TH2}(2). Thus $(2)\Rightarrow (3)$. 

Hence Theorem \ref{TH3} holds.

\begin{corollary}\label{C5}
Assume $R$ is a B\'{e}zout domain such that for each $(a,b)\in Um(R^2)$ and $c\in R$, the image of the functorial homomorphism 
$$U(R/Rabc)\rightarrow U(R/Rab)\cong U(R/Ra)\times U(R/Rb)$$ 
is the product of the images of the following two functorial homomorphisms $\break U(R/Rabc)\rightarrow U(R/Ra)$ and $U(R/Rabc)\rightarrow U(R/Rb)$. Then $R$ is an \textsl{EDD}.
\end{corollary}

\begin{proof}
Based on (the implication $(3)\Rightarrow (1)$ of) Theorem \ref{TH3}, it suffices to show that $R$ is a $U_2$ ring. Given a unit $u+Rc\in U(R/Rc)$ and $(a,b)\in Um(R^{2})$, let $a^{\prime
}:=\gcd (a,c)$ and $b^{\prime }:=\gcd (b,c)$. As $\gcd (a,b)=1$, $a^{\prime }b^{\prime }$ divides $c$ and hence there
exists $c^{\prime }\in R$ such that $c=a^{\prime }b^{\prime }c^{\prime }$.
From our hypothesis applied to $(a^{\prime },b^{\prime })\in Um(R^2)$ and $c^{\prime}\in R$ it follows that there exist units $e+Rc,f+Rc\in U(R/Rc)$ such that their
images in $U(R/Ra^{\prime })\times U(R/Rb^{\prime })$ are $(u+Ra^{\prime
},1+Rb^{\prime })$ and $(1+Ra^{\prime },u+Rb^{\prime })$ (respectively).
Then $(u+Ra,1+Rb)\in R/Rab\cong R/Ra\times R/Rb$ and $e+Rc\in R/Rc$ (resp.\ $(1+Ra,u+Rb)\in R/Rab\cong R/Ra\times R/Rb$ and $f+Rc\in R/Rc$) map to the
same element in $R/Ra^{\prime }b^{\prime }\cong R/Ra^{\prime }\times
R/Rb^{\prime }$ and hence, as $abc^{\prime }$ is the least common multiple
of $c$ and $ab$, there exists a pair $(g,h)\in R^2$ such that $g+Rabc^{\prime }$ (resp.\ $h+Rabc^{\prime }$) reduces to both of them. As we have $g+Rbc\in U(R/Rbc)$ and $h+Rac\in U(R/Rac)$, it follows that $(e+Rc)(f+Rc)\in\Im(\upsilon_{a,b,c})$. Hence to show that $\upsilon_{a,b,c}$ is surjective, by replacing $u+Rc$ with $(u+Rc)[(e+Rc)(f+Rc)]^{-1}$ we can assume that $u-1\in
Ra^{\prime }b^{\prime }$, i.e., the images of $u+Rc$ and $1+Rab$ in $Ra^{\prime}b^{\prime}$ are equal, which implies that $u+Rc\in\Im\bigl(U(R/Rabc^{\prime
}\bigr)\rightarrow U\bigl(R/Rc)\bigr)$. As $abc^{\prime }+Rabc\in R/Rabc$ has square $0$, the homomorphism $U(R/Rabc)\rightarrow U(R/Rabc^{\prime })$ is
surjective and thus $u+Rc\in\Im\bigl(U(R/Rabc)\rightarrow U(R/Rc)\bigr)\subset\Im(\upsilon_{a,b,c})$. So $\Im(\upsilon_{a,b,c})$ is surjective. Thus $R$ is a $U_2$ ring.\end{proof}

\begin{corollary}
\label{C4} Let $R$ be a Hermite ring. Then $R$ is an \textsl{EDR}
iff for all homomorphisms $\phi :\mathbb{Z}[x,y,z]\rightarrow R$,
the image of $\mathcal{D}\in Um\bigl(\mathbb{M}_{2}(\mathbb{Z}[x,y,z])\bigr)$
(equivalently, $\mathcal{F}\in Um\bigl(\mathbb{M}_{2}(\mathbb{Z}[x,y,z])\bigr)$, see
Case 1 of Section \ref{S2}) in $Um\bigl(\mathbb{M}_{2}(R)\bigr)$ via $\phi $ is extendable.
\end{corollary}

\begin{proof}
The `only if' part is obvious. To prove the `if' part, based on (the implication $(3)\Rightarrow (1)$ of) Theorem \ref{TH3}, it suffices to prove that $R$ is a $U_{2}$ ring, i.e., for each $(a,c)\in R^2$ and every $(u,c)\in Um(R^2)$, we have $u+Rc\in \Im(\upsilon_{a,1-a,c})$. Let $(s,t)\in R^2$ be such that $cs+ut=1$. As $u+Rcs\in U(R/Rcs)$, it suffices to prove that $u+Rcs\in \Im(\upsilon_{a,1-a,cs})$. Thus by replacing $c$ with $cs$ we can assume that $c=1-ut$. Hence $A=\left[ 
\begin{array}{cc}
a(1-ut) & u \\ 
0 & (1-a)(1-ut)\end{array}\right] \in Um\bigl(\mathbb{M}_{2}(R)\bigr)$ is the image of $\mathcal{D}$ via the
homomorphism $\phi :\mathbb{Z}[x,y,z]\rightarrow R$ that maps $x$, $y$, and $z$ to $a$, $u$, and $t$ (respectively) and so $A$ is extendable. As $\mathcal{D}$
is equivalent to $\mathcal{E}$ (see Case 1 of Section \ref{S2}), it follows that $A $ is equivalent to the matrix $B:=\left[ 
\begin{array}{cc}
a & u \\ 
0 & (1-a)(1-uw)\end{array}\right] \in Um\bigl(\mathbb{M}_{2}(R)\bigr)$, where $w:=2t-ut^{2}$. From \cite{CPV1}, Lem.\ 4.1(3)
it follows that $B$ is extendable and hence its reduction modulo $\det(B)=Ra(1-a)(1-uw)=Ra(1-a)c^{2}$ is non-full (see \cite{CPV1}, Prop.\ 5.1(2)).
We denote by $\bar{\ast}$ the reduction modulo $Ra(1-a)c^{2}$ of $\ast$, where $\ast $ is either $R$ or an element of $R$. As $B$ is non-full, we have a
product decomposition $\bar{u}=\bar{u_{1}}\bar{u}_2$ with $\bar{u}_{1}=u_{1}+Ra(1-a)c^{2}\in \bar{R}$ dividing $\bar{a}$ and $\bar{u}_{2}=u_{2}+Ra(1-a)c^{2}\in\bar R$ dividing $(1-\bar{a})(1-\bar{w}\bar{u})$ and hence
also $1-\bar{a}$; here $(u_1,u_2)\in R^2$. Thus there exists a triple $(d,e,f)\in R^3$ such that $u_{1}u_{2}=u+a(1-a)c^{2}d$, $u_{1}$ divides $a+a(1-a)c^{2}e$ and $u_{2}$
divides $1-a+a(1-a)c^{2}f$. It follows that $\Bigl(\bigl(u_{1},(1-a)c\bigr),(u_{2},ac)\Bigr)\in
Um(R^{2})^2$, and thus we have an identity $u+Rc=u_{1}u_{2}+Rc=\upsilon_{a,1-a,c}\bigl(u_{2}+Rac,u_{1}+R(1-a)c\bigr).$
\end{proof}

\begin{corollary}\label{C5} Let $R$ be a B\'{e}zout domain. Then $R$ is an \textsl{EDD} iff for all triples $(a,u,t)\in R^3$ with $u\neq 0$ there exists $(s,l,z)\in R^3$ such that 
\begin{equation}\label{EQ4}
(1-us-al)^2+l-usl-al^2-(s+t-ust)z=0.
\end{equation}
\end{corollary}

\begin{proof}
The Hermite ring $R$ is an \textsl{EDR} iff for all triples $(a,u,t)\in R^3$, the matrix $\break A=\left[ 
\begin{array}{cc}
a(1-ut) & u \\ 
0 & (1-a)(1-ut)\end{array}\right] \in Um\bigl(\mathbb{M}_{2}(R)\bigr)$ is extendable (see Corollary \ref{C4}), equivalently it is determinant liftable (see Example \ref{EX1}(1)), and hence equivalently (see \cite{CPV2}, Thm.\ 1.2) there exists $(x,y,z,w)\in R^4$ such that 
$$1-a(1-ut)x-uy-(1-a)(1-ut)w=0=xw-yz.$$
If $u=0$, then $A$ is diagonal and hence simply extendable (see \cite{CPV1}, Ex.\ 4.9(3)).
Thus we can assume $u\neq 0$. The equation $1-a(1-ut)x-uy-(1-a)(1-ut)w=0$ can be rewritten as $1-uy=(1-ut)[ax-(1-a)w]$ and working modulo $Ru$ with $u\neq 0$ it follows that its general solution is $ax+(1-a)w=1-us$ and $y=s+t-ust$ with $s\in R$. As $R$ is an integral domain, the general solution of $ax+(1-a)w=1-us$ is $x=1-us+(1-a)l$ and $w=1-us-al$ with $l\in R$ and the equation $xw-yz=0$ becomes $[1-us+(1-a)l](1-us-al)-(s+t-ust)z=0$ and the corollary follows.
\end{proof}

\section{Proof of Theorem \ref{TH4}}\label{S5}

We first prove two general lemmas.

\begin{lemma}\label{L1}
Let $(d,e)\in R^2$ be such that $Rd=Re$. Then there exists $N\in SL_2(R)$ such that $N\left[\begin{array}{cc}
d & 0 \\ 
0 & 0\end{array}\right]=\left[\begin{array}{cc}
e & 0 \\ 
0 & 0\end{array}\right]$ (so $\left[\begin{array}{cc}
d & 0 \\ 
0 & 0\end{array}\right]$ and $\left[\begin{array}{cc}
e & 0 \\ 
0 & 0\end{array}\right]$ are equivalent).
\end{lemma}

\begin{proof}
Let $(u,v)\in R^2$ be such that $(d,e)=(eu,dv)$; thus $1-uv\in\Ann_R(d)$. For $N:=\left[\begin{array}{cc}
v & -1 \\ 
1-uv & u\end{array}\right]\in SL_2(R)$, one computes $N\left[\begin{array}{cc}
d & 0 \\ 
0 & 0\end{array}\right]=\left[\begin{array}{cc}
e & 0 \\ 
0 & 0\end{array}\right]$.
\end{proof}

\begin{lemma}\label{L2}
Assume $N(R)=0$. Let $e\in R$ be such that for $m\in\{2,4\}$ the reduction function $Um(R^m)\to Um\Bigl(\bigl(R/\Ann_R(e)\bigl)^m\Bigr)$ is surjective\footnote{Each element of $Um\Bigl(\bigl(R/\Ann_R(e)\bigl)^m\Bigr)$ is the image of an element of $Um\bigl((R/Rf)^m\bigr)$ for some $f\in\Ann_R(e)$ and hence, if $sr(R)\le m$, of an element in $Um(R^m)$ by \cite{CPV1}, Prop.\ 2.4(1).}. Then all zero determinant matrices in $eUm\bigl(\mathbb M_2(R)\bigr)$ admit diagonal reduction iff $R/\Ann_R(e)$ is a $\Pi_2$ ring. 
\end{lemma}

\begin{proof}
For $*\in R\cup\mathbb M_2(R)$, let $\bar{*}$ be the reduction of $*$ modulo $\Ann_R(e)$. For two matrices $M=\left[\begin{array}{cc}
m_{11} & m_{12} \\ 
m_{21} & m_{22}\end{array}\right]$ and $L=\left[\begin{array}{cc}
l_{11} & l_{12} \\ 
l_{21} & l_{22}\end{array}\right]$ in $\mathbb M_2(R)$, the product
$$\left[\begin{array}{cc}
m_{11} & m_{12} \\ 
m_{21} & m_{22}\end{array}\right]\left[\begin{array}{cc}
e & 0 \\ 
0 & 0\end{array}\right]\left[\begin{array}{cc}
l_{11} & l_{12} \\ 
l_{21} & l_{22}\end{array}\right]=\left[\begin{array}{cc}
em_{11}l_{11} & em_{11}l_{12} \\ 
em_{21}l_{11} & em_{21}l_{12}\end{array}\right]$$
depends only on the first column of $M$ and the first row of $L$, i.e., it depends only on the quadruple $(m_{11},m_{21},l_{11},l_{12})\in R^4$. Given such a quadruple in $R^4$, we have $\bigl((m_{11},m_{21}),(l_{11},l_{12})\bigr)\in Um(R^2)^2$ iff there exists $(m_{12},m_{22},l_{21},l_{22})\in R^4$ such that the resulting matrices $M$ and $L$ are in $SL_2(R)$. 

Let $A=\left[\begin{array}{cc}
a & b \\ 
c & d\end{array}\right]\in Um\bigl(\mathbb M_2(R)\bigr)$ be such that $\det(eA)=0$; so $e^2\det(A)=0$ and as $N(R)=0$ it follows that $e\det(A)=0$. Thus $\det(A)\in\Ann_R(e)$ and hence $\det(\bar A)=0$. The matrix $eA$ admits diagonal reduction iff it is equivalent to $(f,0)$ (see Example \ref{EX5}) with $Re=Rf$ and hence iff it is equivalent to $(e,0)$ by Lemma \ref{L1}. From this and the previous paragraph it follows that $eA$ admits diagonal reduction iff there exists $\bigl((m_{11},m_{21}),(l_{11},l_{12})\bigr)\in Um(R^2)^2$ such that $(ea,eb,ec,ed)=(em_{11}l_{11},em_{11}l_{12},em_{21}l_{11},em_{21}l_{12})$ and, as the map $Um(R^2)\to Um\Bigl(\bigl(R/\Ann_R(e)\bigl)^2\Bigr)$ is surjective, iff $\bar A=\left[\begin{array}{c}
\overline{m_{11}}\\ 
\overline{m_{21}}\end{array}\right]\left[\begin{array}{cc}
\overline{l_{11}} & \overline{l_{12}}\end{array}\right]$ is non-full. Each $C\in Um\Bigl(\mathbb M_2\bigl(R/\Ann_R(e)\bigl)\Bigl)$ is $\bar B$ for some $B\in Um\bigl(\mathbb M_2(R)\bigl)$ by hypotheses; if $\det(C)=0$, then $\det(B)\in\Ann_R(e)$ and hence $\det(eB)=0$. Thus we conclude that all zero determinant matrices in $eUm\bigl(\mathbb M_2(R)\bigr)$ admit diagonal reduction iff all zero determinant matrices in $Um\Bigl(\mathbb M_2\bigl(R/\Ann_R(e)\bigl)\Bigl)$ are non-full (equivalently, are simply extendable by \cite{CPV1}, Prop.\ 5.1(1)) and iff $R/\Ann_R(e)$ is a $\Pi_2$ ring.\end{proof}

We prove Theorem \ref{TH4}. If $R$ is an \textsl{EDR} then all zero determinant matrices in $\mathbb M_2(R)$ admit diagonal reduction. If for each $a\in R$, $R/Ann_R(a)$ is a $\Pi_2$ ring, then all zero determinant matrices in $RUm\bigl(\mathbb M_2(R)\bigr)=\mathbb M_2(R)$ admit diagonal reduction by Lemma \ref{L2}. As $N(R)=0$, in the last two sentences we can replace `admit diagonal reduction' with `are non-full' by Example \ref{EX5}. Thus to prove Theorem \ref{TH4} we can assume that the Hermite $R$ is such that all zero determinant matrices in $\mathbb M_2(R)$ are non-full, in which case Theorem \ref{TH4} follows from (the equivalence $(1)\Leftrightarrow (6)$ of) Theorem  \ref{TH1}.

\section{On ${(WSU')}_{n}$ and ${(WSU)}_{n}$ rings}\label{S6}

For $F\in\mathbb M_n(R)$, let $F^T\in\mathbb M_n(R)$ be its transpose and let $\Ker_F$ and $\Im_F$ be the kernel and the image (respectively) of the $R$-linear map $L_F:R^n\rightarrow R^n$ defined by it. Note that, as usual, the elements of the domain or codomain of $L_F$ are written as column vectors.

\begin{proposition}\label{P4}
A ring $R$ is an ${(SU^{\prime})}_{n}$ (resp.\ ${(SU)}_{n}$) ring in the sense of Definition \ref{def2}(2) iff it is so in the sense of \cite{lor}, Def.\ 2.1. 
\end{proposition}

\begin{proof}
For the `only if' part, let $B\in\mathbb M_n(R)$. As $R$ is a Hermite ring, there exist $e\in R$ and $C\in Um\bigl(\mathbb M_n(R)\bigr)$ such that $B=eC$. As $R$ is an ${(SU^{\prime})}_{n}$ (resp.\ ${(SU)}_{n}$) ring, there exists $N\in SL_n(R)$ (resp.\ $N\in GL_n(R))$ such that $CN$ is symmetric. Hence $BN=e(CN)$ is symmetric, thus $R$ is an ${(SU^{\prime})}_{n}$ (resp.\ ${(SU)}_{n}$) ring in the sense of \cite{lor}, Def.\ 2.1. 

For the `if' part, clearly $R$ is a ${(WSU^{\prime})}_{n}$ (resp.\ ${(WSU)}_{n}$) ring and it is a Hermite ring by \cite{lor}, Prop.\ 3.1. 
\end{proof}

Directly from Proposition \ref{P4} and \cite{lor}, Prop.\ 3.3 we get the following result. 

\begin{corollary}\label{C5.5}
Let $n\in\mathbb N\setminus\{1,2\}$. Assume $R$ is an integral domain and an ${(SU)}_{n}$ ring. Then $R$ is an ${(SU)}_m$ ring for each $m\in\{2,\ldots,n-1\}$.
\end{corollary}

As in \cite{lor}, Sect.\ 2, if $U(R)$ is a $2$-divisible abelian group, $R$ is a ${(WSU^{\prime})}_{n}$ ring iff it is a ${(WSU)}_{n}$ ring. For $(WSU)_n$ rings we have the following general result.

\begin{theorem}\label{TH7}
Assume $R$ is a ${(WSU)}_{n}$ ring. Then the following properties hold.

\medskip \textbf{(1)} Each unimodular matrix $A\in Um\bigl(\mathbb M_n(R)\bigr)$ is equivalent to its transpose $A^T$. 

\smallskip \textbf{(2)} If $n=2$, then every projective $R$-module of rank $1$ and generated by $2$ elements is self-dual.

\smallskip \textbf{(3)} If $sr(R)\le n^2$, then for each $d\in R$, $R/Rd$ is a ${(WSU)}_n$ ring.\end{theorem}

\begin{proof}
If $N\in GL_n(R)$ is such that $AN$ is symmetric, i.e., $AN=N^TA^T$ then $(N^T)^{-1}AN=A^T$ and hence part (1) holds. 

To check part (2), let $P$ be a projective $R$-module of rank $1$ generated by $2$ elements. We consider an isomorphism $P\oplus Q\cong R^2$ to be viewed as an identification. Taking determinants it follows that $Q$ is the dual of $P$. Let $B\in Um\bigl(\mathbb M_2(R)\bigr)$ be such that $\Ker_{B}=Q$ and $\Im_{B}=P$. As $B^T$ and $B$ are equivalent by part (1), the $R$-modules $P$ and $Q$ are isomorphic, so $P$ is self-dual. Thus part (2) holds. 

To check part (3) for $d\in R$, let $\bar A\in Um\bigl(\mathbb M_n(R/Rd)\bigr)$. As $sr(R)\le n^2$, there exists $A\in Um\bigl(\mathbb M_n(R)\bigr)$ such that its reduction modulo $Rd$ is $\bar A$ (see \cite{CPV1}, Prop.\ 2.4(1)). Let $N\in GL_n(R)$ be such that $AN$ is symmetric. If $\bar N\in GL_n(R/Rd)$ is the reduction of $N$ modulo $Rd$, then $\bar A\bar N$, being the reduction of $AN$ modulo $Rd$, is symmetric. Thus $R/Rd$ is a ${(WSU)}_n$ ring; so part (3) holds.
\end{proof}

\begin{remark}\normalfont\label{rem1} Assume $R$ is an integral domain such that for each $d\in R$, every projective $R/Rd$-module of rank $1$ and generated by $2$ elements is self-dual and the functorial homomorphism $GL_2(R)\rightarrow GL_2(R/Rd)$ is surjective. Then each matrix $A\in Um\bigl(\mathbb M_2(R)\bigr)$ is equivalent to its transpose $A^T$.
\end{remark}

\begin{definition}\label{def4} Let $(m,n)\in\mathbb N^2$ be such that $m<n$. We say that $R$ is:

\medskip
{\bf (1)} a $WH_{n,m}^{\prime}$ (resp.\ $WH_{n,m}$) ring if for each $(A_1,\ldots,A_m)\in Um\bigl(\mathbb M_n(R)\bigr)^m$ there exists $N\in SL_n(R)$ (resp.\ $N\in GL_n(R)$) such that the trace of $A_iN$ is $0$ for every $i\in\{1,\ldots,m\}$.

\smallskip
{\bf (2)} an $SH_{n,m}^{\prime}$ (resp.\ $SH_{n,m}$) ring if it is a Hermite ring and a $WH_{n,m}^{\prime}$ (resp.\ $WH_{n,m}$) ring.

\smallskip
{\bf (3)} an $H_{n,m}^{\prime}$ (resp.\ $H_{n,m}$) ring if for each $(A_1,\ldots,A_m)\in\mathbb M_n(R)^m$ there exists $N\in SL_n(R)$ (resp.\ $N\in GL_n(R)$) such that the trace of $A_iN$ is $0$ for every $i\in\{1,\ldots,m\}$.
\end{definition}

Definition \ref{def4}(3) was introduced by Lorenzini im \cite{lor}, Def.\ 2.4.  In Definition \ref{def4}(1) or (3), the existence of $N$ is equivalent to the existence of $M\in {SL}_n(R)$ (resp.\ $M\in {GL}_n(R)$) such that the trace of $MA_i$ is
$0$ for every $i\in\{1,\ldots,m\}$ and equivalent to the existence of $(M,N)\in {SL}_n(R)^2$ (resp.\ $(M,N)\in {GL}_n(R)^2$) such that the trace of $MA_iN$ is
$0$ for every $i\in\{1,\ldots,m\}$. This follows via conjugation from the fact that traces are preserved under transposition.

Clearly, if $n-1\ge m\ge 2$ and $l\in\{1,\ldots,m-1\}$, then each $WH_{n,m}^{\prime}$ (or $WH_{n,m}$ or $SH_{n,m}^{\prime}$ or $SH_{n,m}$) ring is a $WH_{n,l}^{\prime}$ (or $WH_{n,l}$ or $SH_{n,l}^{\prime}$ or $SH_{n,l}$) ring.

The trace of the product matrix $\left[\begin{array}{cc}
a & b \\ 
c & d\end{array}\right]\left[\begin{array}{cc}
0 & -1 \\ 
1 & 0\end{array}\right]$ is $b-c$. Based on this it follows easily that $R$ is a $WH_{2,1}^{\prime}$ (resp.\ $WH_{2,1}$) ring iff it is a ${(WSU^{\prime})}_2$ (resp.\ ${(WSU)}_2$) ring (cf.\ \cite{lor}, Prop.\ 4.11).

\begin{proposition}\label{P4.5}
Let $n\in\mathbb N\setminus\{1\}$. Then the following properties hold.

\medskip
{\bf (1)} Let $m\in\{1,\ldots,n-1\}$. If $R$ is an $SH_{n,m}^{\prime}$ (resp.\ $SH_{n,m}$) ring, then it is an $H_{n,m}^{\prime}$ (resp.\ $H_{n,m}$) ring. 

\smallskip
{\bf (2)} A ring $R$ is an $SH_{n,n-1}^{\prime}$ (resp.\ $SH_{n,n-1}$) ring iff it is an $H_{n,n-1}^{\prime}$ (resp.\ $H_{n,n-1}$) ring.

\smallskip
{\bf (3)} A ring $R$ is an $SH_{2,1}^{\prime}$ ring iff it is an $H_{2,1}^{\prime}$ ring and iff it is an ${(SU^{\prime})}_2$ ring.

\smallskip
{\bf (4)}  A ring $R$ is an $SH_{2,1}$ ring iff it is an $H_{2,1}$ ring and iff it is an ${(SU)}_2$ ring.
\end{proposition}

\begin{proof}
For part (1), let $(B_1,\ldots,B_m)\in\mathbb M_n(R)^m$. As $R$ is a Hermite ring, for each $i\in\{1,\ldots,m\}$ there exists $(e_i,C_i)\in R\times Um\bigl(\mathbb M_n(R)\bigr)$ such that $B_i=e_iC_i$. As $R$ is an $WH_{n,m}^{\prime}$ (resp.\ $WH_{n,m}$) ring, there exists $N\in SL_n(R)$ (resp.\ $N\in GL_n(R))$ such that the trace of $C_iN$ is $0$ for every $i\in\{1,\ldots,m\}$. Hence the trace of $B_iN=e_i(C_iN)$ is $0$ for every $i\in\{1,\ldots,m\}$. Thus $R$ is an $H_{n,m}$ (resp.\ $H_{n,m}^{\prime}$) ring. So part (1) holds.

Based on part (1), to prove part (2) it suffices to prove its `if' part. As $R$ is an $H_{n,n-1}^{\prime}$ (resp.\ $H_{n,n-1}$) ring, it is also a $WH_{n,n-1}^{\prime}$ (resp.\ $WH_{n,n-1}$) and it is a Hermite ring by \cite{lor}, Lem.\ 4.3. So $R$ is an $SH_{n,n-1}^{\prime}$ (resp.\ $SH_{n,n-1}$) ring. Thus part (2) holds.

Parts (3) and (4) follow from part (2) applied to $n=2$ and the paragraph before this proposition.
\end{proof}

\begin{remark}\normalfont\label{rem1.5} {\bf (1)} For $(a,b,c,d)\in Um(R^4)$, let 
$$\mathcal L_{(a,b,c,d)}:=\{(ax+by+cz+dw,xw-yz)|(x,y,z,w)\in R^4\}\subset R^2.$$
We have $(1,0)\in \mathcal L_{(a,b,c,d)}$ iff the matrix $\left[\begin{array}{cc}
a & b \\ 
c & d\end{array}\right]\in Um\bigl(\mathbb M_2(R)\bigr)$ is determinant liftable. Moreover, $R$ is a $WJ_{2,1}$ ring in the sense of \cite{CPV2}, Def.\ 1.10(1) iff for each $(a,b,c,d)\in Um(R^4)$ we have $\mathcal L_{(a,b,c,d)}=R^2$. Similarly, $R$ is a $WH_{2,1}^{\prime}$ (resp.\ $WH_{2,1}$) ring iff for each $(a,b,c,d)\in Um(R^4)$ we have $(0,1)\in\mathcal L_{(a,b,c,d)}$ (resp.\ $[\{0\}\times U(R)]\cap \mathcal L_{(a,b,c,d)}\neq\emptyset$), cf.\ \cite{lor}, Lem.\ 4.1.
 
\smallskip
{\bf (2)} The sets $\mathcal L_{(a,b,c,d)}$ of part (1) have many properties. For instance, if $u\in U(R)$ and $(\Psi,\Delta)\in \mathcal L_{(a,b,c,d)}$, then $\bigl(u\Psi+2(ad-bc),u^2\Delta+u\Psi+ad-bc\bigr)\in \mathcal L_{(a,b,c,d)}$; this is so as for $(x,y,z,w)\in R^4$ with $(ax+by+cz+dw,xw-yz)=(\Psi,\Delta)$ for $(x^{\prime},y^{\prime},z^{\prime},w^{\prime}):=(ux+d,uy-c,uz-b,uw+a)$ we have 
$$(ax^{\prime}+by^{\prime}+cz^{\prime}+dw^{\prime},x^{\prime}w^{\prime}-y^{\prime}z^{\prime})=\bigl(u\Psi+2(ad-bc),u^2\Delta+u\Psi+ad-bc\bigr).$$ A similar argument gives that if $\bigl((1,t),(0,u)\bigr)\in\mathcal L_{(a,b,c,d)}^2$, then there exists $s\in R$ such that for each $x\in R$ we have $(1,ux^2+sx+t)\in\mathcal L_{(a,b,c,d)}^2$.
\end{remark}

\section{Proof of Theorem \protect\ref{TH5} and applications}\label{S7}

Let $\Spec R$ be the spectrum of $R$ and let $\Max R$ be its subset of maximal ideals.

To prove part (1), it suffices to show that for all $(u,c)\in Um(R^2)$ and $a\in R$, we
have $u^{2}+Rc\in\Im(\upsilon_{a,1-a,c})$. To check this we can assume that $c\neq
0 $. As in the proof of Corollary \ref{C4}, we can assume that there exists $t\in R$
such that $c=1-ut$. As we are assuming that $R$ is a $({WSU})_{2}$ ring, there
exists $M=\left[ 
\begin{array}{cc}
x & y \\ 
z & w\end{array}\right] \in {GL}_{2}(R)$ such that $M\left[ 
\begin{array}{cc}
ac & u \\ 
0 & (1-a)c\end{array}\right] $ is a symmetric matrix, i.e., we have an identity
\begin{equation}
acz=ux+(1-a)cy. \label{EQ6}
\end{equation}As $c$ divides $ux$ and $(u,c)\in Um(R^{2})$, $c$ divides $x$. Let $s\in R$ be such that $x=cs$.

Let $\mu:=\det (M)=csw-yz\in U(R)$. Thus $(cs,yz)\in Um(R^{2})$. 

We first assume that $c\notin Z(R)$. Dividing Equation
(\ref{EQ6}) by $c$ it follows that 
\begin{equation}
az+(a-1)y=us. \label{EQ7}
\end{equation}

We check that the assumption $(a,s)\notin Um(R^2)$ leads to a contradiction. This assumption implies that there exists $\mathfrak m\in\Max R$ such that $Ra+Rs\subseteq\mathfrak m$. From this and Equation (\ref{EQ7}) it follows that $(a-1)y\in\mathfrak m$. As $a-1\not\in\mathfrak m$, it follows that $y\in\mathfrak m$, hence $(s,yz)\notin Um(R^2)$, a contradiction. As $(a,s)\in Um(R^2)$, $s+Ra\in U(R/Ra)$.

From Equation (\ref{EQ7}) it follows that the images of $(-y+Rc)(u+Rc)^{-1}\in U(R/Rc)$ and $s+Ra\in U(R/Ra)$ in $U\bigl(R/(Ra+Rc)\bigr)$ are equal. Thus there exists a unit of the quotient ring $R/(Ra\cap Rc)$ that maps to 
$(-y+Rc)(u+Rc)^{-1}\in U(R/Rc)$. 
As the ideal $(Ra\cap Rc)/Rac$ of $R/Rac$ has square $0$, the functorial homomorphism $U(R/Rac)\rightarrow U\bigl(R/(Ra\cap Rc)\bigr)$ is surjective. From the last two sentences it follows that $(-y+Rc)(u+Rc)^{-1}\in \Im(\upsilon_{a,1-a,c})$.

Due to the symmetry of Equation (\ref{EQ7}) in $az$ and $(a-1)y$, similar arguments give $s+R(1-a)\in U\bigl(R/R(1-a)\bigr)$ and $(z+Rc)(u+Rc)^{-1}\in \Im(\upsilon_{a,1-a,c})$. Clearly $\mu+Rc\in \Im(\upsilon_{a,1-a,c})$. As $\mu+Rc=(-y+Rc)(z+Rc)$, it follows that 
\begin{equation*}
(u+Rc)^{2}=[(-y+Rc)(u+Rc)^{-1}]^{-1}[(z+Rc)(u+Rc)^{-1}]^{-1}(\mu+Rc)\in \Im(\upsilon_{a,1-a,c}).
\end{equation*}

We use a trick with annihilators to show that an analogue of Equation (\ref{EQ7}) holds even if $c\in Z(R)$. As $\bigl((c,u),(c,yz)\bigr)\in Um(R^2)^2$, it follows that $(c,uyz)\in Um(R^2)$. This implies
that $\Ann_R(c)\subset Ruyz$. As $az+(a-1)y-us\in\Ann_R(c)$ by Equation (\ref{EQ6}), it follows that there exists $q\in R$ such that 
$az+(a-1)y-us=uyzq$. If $s_1:=s+yzq$, then $(s_1,yz)\in
Um(R^2)$ and the identity $az+(a-1)y=us_1$ holds. Thus, using $s_1$ and the
last identity instead of $s$ and Equation (\ref{EQ7}), as above we
argue that $(u+Rc)^{2}\in \Im(\upsilon_{a,1-a,c})$. So part (1) holds.

Part (2) follows from Theorem \ref{TH7}(2) and (3). 

The first part of part (3) follows from part (1) and definitions. Based on it, the second part of part (3) follows from (the implication $(3)
\Rightarrow (1)$ of) Theorem \ref{TH3}. 

Thus Theorem \ref{TH5} holds.

\begin{remark}\normalfont\label{rem2}
If we have $(a,y)\in Um(R^2)$ (resp.\ $(1-a,z)\in Um(R^2)$), then $-y+Rac\in U(R/Rac)$ (resp.\ $z+R/R(1-a)c\in U\bigl(R/R(1-a)c\bigr)$), and therefore $u+Rc\in \Im(\upsilon_{a,1-a,c})$.
\end{remark}

\begin{corollary}
\label{C6} For a B\'{e}zout domain the following statements are equivalent.

\medskip \textbf{(1)} The ring $R$ is a $(WSU)_2$ (resp.\ $(WSU^{\prime})_2$) ring.

\smallskip \textbf{(2)} For all $\bigl((a,b),(c,u)\bigr)\in Um(R^{2})^2$ with $abc\neq 0$,
if a pair $(a^{\prime},b^{\prime})\in R^2$ such that $aa^{\prime}+bb^{\prime}=1$ is given, then there exists a triple $(s,l,w)\in R^3$ with the property that $csw+(a^{\prime}us-bl)(b^{\prime}us+al)\in U(R)$ (resp.\ $csw+(a^{\prime}us-bl)(b^{\prime}us+al)=1$).\end{corollary}

\begin{proof}
If $A_1,A_2\in\mathbb M_2(R)$ are equivalent, then there exists $M_1\in {GL}_2(R)$ such that $M_1A_1$ is symmetric iff there exists $M_2\in {GL}_2(R)$ such that $M_2A_2$ is symmetric. This is so, as for $(M,N,M_1)\in {GL}_2(R)^3$ such that $A_2=MA_1N$ and $M_1A_1$ is symmetric, by denoting $M_2:=N^TM_1M^{-1}\in {GL}_2(R)$, $M_2A_2=N^T(M_1A_1)N$ is symmetric. Also, if $(M,N,M_1)\in {SL}_2(R)^3$, then $M_2\in {SL}_2(R)$. 

Based on the previous paragraph and Case 1 of Section \ref{S2}, the ring $R$ is a $(WSU)_2$ (resp.\ $(WSU^{\prime})_2$) ring iff for all $\bigl((a,b),(c,u)\bigr)\in Um(R^{2})^2$ there exists $M=\left[ 
\begin{array}{cc}
x & y \\ 
z & w\end{array}\right]$ in ${GL}_2(R)$ (resp.\ in ${SL}_2(R)$) such that the matrix $M\left[ 
\begin{array}{cc}
ac & u \\ 
0 & bc\end{array}\right] $ is symmetric, i.e., the equation $acz=ux+bcy$ holds. If $abc=0$, then $\left[ 
\begin{array}{cc}
ac & u \\ 
0 & bc\end{array}\right] $ is simply extendable (see \cite{CPV1}, Ex.\ 4.9(3)) and hence equivalent to a diagonal matrix (see \cite{CPV1}, Thm.\ 4.3). From this and the previous paragraph we get that $M$ exists. Thus we can assume that $abc\neq 0$. So $c\notin Z(R)=\{0\}$ and hence as in the proof of Theorem \ref{TH5}, we write $x=cs$ with $s\in R$ and the identity $az=us+by$ holds. The general solution of the equation $az-by=us$ in the indeterminates $y,z$ is $z=a^{\prime}us-bl$ and $y=-b^{\prime}us-al$
with $l\in R$. As $\det(M)=xw-yz=csw+(a^{\prime}us-bl)(b^{\prime}us+al)$, the corollary holds.
\end{proof}

\begin{remark}\normalfont\label{rem3}
Referring to Corollary \ref{C6}(2) with $b=1-a$, we can take $a^{\prime}=1=b^{\prime}$, and one is searching for $(s,l,w)\in R^3$
such that $csw+(us+al-l)(us+al)\in U(R)$ (resp.\ and $csw+(us+al-l)(us+al)=1$,
 which for $c=1-ut$ with $t\in R$ is closely related to Equation (\ref{EQ4}) via a substitution of the form $z:=-w$).
\end{remark}

\begin{remark}\normalfont\label{rem4}
A symmetric matrix $A=\left[\begin{array}{cc}
a & b \\ 
b & d\end{array}\right]\in Um\bigl(\mathbb M_2(R)\bigr)$ is determinant liftable iff there exists $(x,w)\in R^2$ such that $b$ divides $1-ax-dw$ and for a suitable $s\in R$ such that $bs=1-ax-cw$, the equation $X^2+sX+xw=0$ has two solutions $y,z\in R$ (so $xw-yz=0$ and $1-ax-by-bz-dw=0$).
\end{remark}

Recall from \cite{ZR}, Sect.\ 1 that $R$ is said to have square stable range $1$, and one writes $ssr(R)=1$, if for each $(a,b)\in Um(R^{2})$, there exists $r\in R$ such that $a^{2}+br\in U(R)$. Clearly, we have $ssr(R)=1$ iff for each $b\in R$, the cokernel of the reduction homomorphism 
$U(R)\rightarrow U(R/Rb)$ is a Boolean group. From this and \cite{CPV1}, Prop.\ 2.4(1) it follows that if $sr(R)=1$, then $ssr(R)=1$. 

\begin{proposition}\label{EX11}
Assume $ssr(R)=1$. Then the following properties hold.

\medskip
{\bf (1)} For $\bigl((a,b),(c,u)\bigr)\in Um(R^2)^2$, let $A:=\left[ 
\begin{array}{cc}
ac & u \\ 
0 & (1-a)c\end{array}\right]$. Then there exists $M\in GL_2(M)$ such that $MA$ is symmetric. 

\smallskip
{\bf (2)} Assume also that $R$ is a Hermite ring. Then each matrix in $Um\bigl(\mathbb{M}_{2}(R)\bigr)$ has a companion test matrix equivalent to a symmetric matrix in $Um\bigl(\mathbb{M}_{2}(R)\bigr)$. Moreover, $R$ is an \textsl{EDR} iff each symmetric matrix in $Um\bigl(\mathbb M_2(R)\bigr)$ is extendable.
\end{proposition}

\begin{proof}
For part (1), as $ssr(R)=1$, there exists $w\in R$ be such that $u^{2}+cw\in U(R)$. For $M:=\left[ 
\begin{array}{cc}
c & -u \\ 
u & w\end{array}\right] \in {GL}_{2}(R)$, the product $MA\in Um\bigl(\mathbb M_2(R)\bigr)$ is symmetric. So part (1) holds.

For part (2), each matrix in $Um\bigl(\mathbb{M}_{2}(R)\bigr)$ has a companion test matrix of the form $A$ of part (1) by Case 1 of Section \ref{S2} and hence which is equivalent to a symmetric matrix in $Um\bigl(\mathbb{M}_{2}(R)\bigr)$ by part (1). Based on this, (the equivalence $(1)\Leftrightarrow (3)$ of) Theorem \ref{TH1} and Proposition \ref{P1} applied to $\mathcal P$ being extendable, it follows that the last sentence of part (2) holds as well.
\end{proof}

\section{Proof of Theorem \ref{TH6}}\label{S8}

Let $P$ be a projective $R/Rd$-module of rank $1$ generated by $2$ elements. We consider an isomorphism $P\oplus Q\cong (R/Rd)^2$. Let $D\in Um\bigl(\mathbb M_2(R/Rd)\bigr)$ be such that $\Ker_D=Q$ and $\Im_D=P$. As $sr(R)\le 2$, there exists $A\in Um\bigl(\mathbb M_2(R)\bigr)$ such that its reduction modulo $Rd$ is $D$ (see \cite{CPV1}, Prop.\ 2.4(1)).
Let $\bar{R}:=R/R\det (A)$ and let $\bar{A}\in Um\bigl(\mathbb{M}_{2}(\bar{R})\bigr)$
be the reduction of $A$ modulo $R\det(A)$. We have $\det(A)\in Rd$. We know that $\Ker_{\bar{A}}$ and $\Im_{\bar{A}}$
are dual projective $\bar{R}$-modules of rank $1$ generated by $2$ elements (see \cite{CPV1}, Lem.\ 3.1) and their reductions modulo the ideal $Rd/R\det(A)$ of $\bar R$ are $Q$ and $P$. Thus, to complete the proof of Theorem \ref{TH6}, it suffices to show that the classes of $\Ker_{\bar{A}}$ and $\Im_{\bar{A}}$ in $\Pic(\bar{R})$ are equal (equivalently, have orders $1$ or $2$ or equivalently, $\Ker_{\bar{A}}^{2}$ and $\Im_{\bar{A}}^{2}$ are free $R$-modules
of rank $2$). So, $d$ and $D$ will not be mentioned again. 

Based on Case 1 of Section \ref{S2}, we can assume that $A=\left[ 
\begin{array}{cc}
ac & u \\ 
0 & bc\end{array}\right] \in \mathbb{M}_{2}(R)$ with $(a,b,c)\in R^3$ such that $\bigl((a,b),(c,u)\bigr)\in Um(R^{2})^2$ and we only use that $\Coker(\upsilon_{a,b,c})$ is a Boolean group. As $\Coker(\upsilon_{a,b,c})$ makes sense for $R/Rabc$ and as two finitely generated projective $R/Rabc^2$-modules are isomorphic if and only if their reductions modulo the square $0$ ideal $Rabc/Rabc^2$ are isomorphic, by reducing modulo $Rabc$ we can assume that $abc=0$. Hence $(R,A)=(\bar R,\bar A)$ and we only use $(R,A)$. 

Let $(s,t)\in
R^{2}$ be such that $cs+ut=1$. Recall that $\Ker_{A}=Rv_{1}+Rv_2$, where $v_1:=\left[ 
\begin{array}{c}
-u \\ 
ac\end{array}\right] $ and $v_2:=\left[ 
\begin{array}{c}
-bc \\ 
0\end{array}\right] $, by \cite{CPV1}, Ex.\ 5.2. As $(a,b)\in Um(R^{2})$, we have $\Spec R=\Spec R_{a}\cup\Spec R_{b}$. Note $(\Ker_{A})_{a}=R_{a}v_{1}$ and $(\Ker_{A})_{b}=R_{b}\left(tv_{1}+\dfrac{s}{b}v_2\right)$. Over $R_{ab}$ we have $c=0$ and hence $v_{1}=u\left(tv_{1}+\dfrac{s}{b}v_2\right)$. Thus we have identities $(\Ker_{A}^{\otimes
2})_a=R_a(v_{1}\otimes v_{1})$, $(\Ker_{A}^{\otimes 2})_b=R_b\left[\left(tv_{1}+\dfrac{s}{b}v_2\right)\otimes \left(tv_{1}+\dfrac{s}{b}v_2\right)\right]$, and over $R_{ab}$ we have an identity $v_{1}\otimes v_{1}=u^{2}\left[\left(tv_{1}+\dfrac{s}{b}v_2\right)\otimes \left(tv_{1}+\dfrac{s}{b} v_2\right)\right]$.
As $u^2+Rc\in\Im(\upsilon_{a,b,c})$, there exists $u_{a}\in U(R/Rbc)$ and $u_{b}\in U(R/Rac)$ such that the product
of their images in $U(R/Rc)$ is $u^{2}+Rc$. As $R_{a}$ is a localization of $R/Rbc$ and $R_{b}$ is a localization of $R/Rac$, we denote also by $u_{a}$ and $u_{b}$ their images in $R_{a}$ and $R_{b}$ (respectively). Thus $u_{a}^{-1}(v_{1}\otimes v_{1})$ and $u_{b}^{-1}\left[\left(tv_{1}+\dfrac{s}{b}v_2\right)\otimes \left(tv_{1}+\dfrac{s}{b}v_2\right)\right]$ coincide over $R_{ab}$ and hence
there exists an element $v\in \Ker_{A}^{\otimes 2}$ whose images in $(\Ker_{A}^{\otimes 2})_{a}$ and $(\Ker_{A}^{\otimes 2})_{b}$ are equal to $u_{a}^{-1}(v_{1}\otimes v_{1})$ and $u_{b}^{-1}\left[\left(tv_{1}+\dfrac{s}{b}v_2\right)\otimes \left(tv_{1}+\dfrac{s}{b}v_2\right)\right]$ (respectively). Thus $Rv=\Ker_{A}^{\otimes 2}$, hence $\Ker_{A}^{\otimes 2}\cong R$.

So Theorem \ref{TH6} holds.

\begin{corollary}
\label{C7} Let $R$ be a Hermite ring such that for each $a\in R$, every self-dual projective $R/Ra$-module of rank $1$ generated by $2$ elements is free. Then the following properties hold.

\medskip \textbf{(1)} The ring $R$ is an \textsl{EDR} iff for all $(a,b)\in Um(R^2)$ and $c\in R$, $\Coker(\upsilon_{a,b,c})$ is a Boolean group. 

\smallskip \textbf{(2)} If $R$ is a $(WSU)_2$ ring, then $R$ is an \textsl{EDR}. 
\end{corollary}

\begin{proof}
The `only if' of part (1) follows from (the implication $(1)\Rightarrow (3)$ of) Theorem \ref{TH3}. To check the `if' of part (1), from Theorem \ref{TH6} and hypotheses it follows that for each $A\in Um\bigl(\mathbb M_2(R)\bigr)$, the reduction $\bar{A}$ of $A$ modulo $R\det(A)$ is such that $\Im_{\bar A}\cong R/R\det(A)$, hence $\bar A$ is non-full and thus extendable (see \cite{CPV1}, Prop.\ 5.1(2)). Hence $R$ is an $E_2$ ring and thus an \textsl{EDR} by (the implication $(3)\Rightarrow (1)$ of) Theorem \ref{TH1}. Thus part (1) holds. 

Part (2) follows from part (1) and Theorem \ref{TH5}(1) (it also follows from Theorems \ref{TH7}(2) and \ref{TH1}).\end{proof}

\section{Proof of Criteria \ref{CR1} and \ref{CR2}}\label{S9}

We first prove Criterion \ref{CR1}(1). If there exists a pair $(e,f)\in R^2$ such that $ae^{2}-cf^{2}\in U(R)$, then from the identity $ae^{2}-cf^{2}=e(ae+bf)-f(be+cf)$ it follows that $(ae+bf,be+cf)\in Um(R^2)$, hence $A$ is simply extendable by \cite{CPV1}, Thm.\ 4.3. For the converse, assume $R$ has characteristic $2$ and $A$ is simply extendable. Let $(e,f)\in Um(R^2)$ be such that $(ae+bf,be+cf)\in Um(R^{2})$ by loc.\ cit. To prove that $ae^{2}-cf^{2}\in U(R)$ it suffices to show that the assumption that there
exist $\mathfrak{m}\in\Max R$ such that $ae^{2}-cf^{2}\in \mathfrak{m}$, leads to a contradiction. By replacing $R$ with $R/\mathfrak{m}$, we can
assume that $R$ is a field and we know that either $a\neq 0$ or $c\neq 0$, either $ae+bf\neq 0$ or $be+cf\neq 0$, $ac=b^{2}$ and $ae^{2}=cf^{2}$;
hence either $e\neq 0$ or $f\neq 0$. If $acef=0$, then by the symmetry in 
$(a,c)$ and $(e,f)$ we can assume that $ae=0$; if $a=0$, then $b=0$, $c\neq 0
$, and hence $f=0$, thus $ae+bf=be+cf=0$, a contradiction,
and if $e=0$, then $f\neq 0$, hence $c=0$, and by the symmetry in the pair $(a,c)$, we similarly reach a contradiction. Thus we can assume that $acef\neq
0$, hence $b^{2}=ac=e^{-2}c^2f^{2}=f^{-2}a^2e^2$; as $R$ is a field of characteristic $2$, it follows that $b=e^{-1}cf=f^{-1}ae$, thus 
$ae+bf=be+cf=0$, a contradiction. So part (1) holds.

To prove Criterion \ref{CR1}(2), as $b\notin Z(R)$, for each $\mathfrak{m}\in\Max R$, $b$
is a nonzero element of $R_{\mathfrak{m}}$. As $R$ is a Hermite ring and $N(R)=0$, each $R_{\mathfrak{m}}$ is a valuation domain (e.g., see footnote of \cite{CPV2}, Thm.\ 11.3(3) for a proof). Let $d\in R$ and $(a^{\prime},b^{\prime})\in Um(R^2)$ be such that $a=da^{\prime}$ 
and $b=db^{\prime}$. As $d$ divides $a$ and $(a,c)\in Um(R^2)$, 
we have $(d,c)\in Um(R^2)$. Thus also $(d^2,c)\in Um(R^2)$. From this, 
as $d^2$ divides $ac$, it follows that $d^2$ divides $a$; let $u\in R$ be such that $a=d^2u$. As $b\in R_{\mathfrak{m}}\setminus\{0\}$, it
is easy to see that the image of $u$ in $R_{\mathfrak{m}}$ is a unit 
for all $\mathfrak{m}\in\Max R$. Hence $u\in U(R)$. By symmetry, there exist $g\in R$ and $v\in U(R)$
such that $c=g^2v$.  As $b\notin Z(R)$ and $b^2=ac=d^2g^2uv$, we have $d\notin Z(R)$ and $g\notin Z(R)$.

To complete the proof that $A$ is extendable, based on Criterion \ref{CR1}(1)
it suffices to show that there exists $(e,f)\in R^{2}$ such that $ae^{2}-cf^{2}=ae^{2}+cf^{2}\in U(R)$. Replacing $A$ by $u^{-1}A$ (see
\cite{CPV1}, Lem.\ 4.1(3)), we can assume that $u=1$. Thus $a=d^{2}$. From this and the
identities $ac=b^{2}$ and $c=g^{2}v$, the element $w:=\dfrac{b}{gd}$ of the total ring of fractions of $R$ has the property that $w^2=v$ in each valuation domain $R_{\mathfrak{p}}$ with $\mathfrak p\in\Spec R$, hence $w\in R_{\mathfrak{p}}$, and we conclude that $w\in R$ and $v=w^2$. By replacing $(v,g)$ with $(1,gw)$, 
we can assume that $v=1$ and $c=g^{2}$. As $(a,c)=(d^{2},g^{2})\in Um(R^{2})$, it follows that $(d,g)\in Um(R^{2})$,
hence there exists $(e,f)\in Um(R^{2})$ such that $de+gf=1$, and so $ae^{2}+cf^{2}=1\in U(R)$. Thus Criterion \ref{CR1}(2) holds. This completes the proof of Criterion \ref{CR1}.

We prove Criterion \ref{CR2}. As $R$ is an $(SU)_{2}$ ring, it is a Hermite ring, and, as the extendability is preserved under matrix equivalence (see \cite{CPV1}, Lem.\ 4.1(3)), it is an $E_{2}$ ring iff each symmetric matrix in $Um\bigl(\mathbb M_2(R)\bigr)$ is extendable. Thus $R$ is an 
\textsl{EDR} iff each symmetric matrix in $Um\bigl(\mathbb M_2(R)\bigr)$ is extendable by Theorem \ref{TH1}. Let $A=\left[ 
\begin{array}{cc}
a & b \\ 
b & c\end{array}\right] \in Um\bigl(\mathbb{M}_{2}(R)\bigr)$. Based on \cite{CPV1}, Lem.\ 4.1(1) and Criterion \ref{CR1}(1) applied to $A$ modulo $R(ac-b^2)$, $A$ is extendable if there exists $(e,f)\in R^{2}$ such that $ae^{2}-cf^{2}+R(ac-b^{2})\in U\bigl(R/[R(ac-b^{2})]\bigr)$, i.e., $(ae^{2}-cf^{2},ac-b^{2})\in Um(R^{2})$, and the converse holds if $R$ has characteristic $2$; hence Criterion \ref{CR2} holds.

\section{Proof of Criterion \protect\ref{CR3}}\label{S10}

As $R$ is a B\'{e}zout domain, $R$ is an \textsl{EDD} iff for all $(a,b,s)\in R^3$, the equation
\begin{equation}
\lbrack 1-(1-a-bs)w-by](1-ax)=ay[(1-a-bs)z+bx] \label{EQ8}
\end{equation}in the indeterminates $x$, $y$, $z$ and $w$ has a solution in $R^{4}$ (see Corollary \ref{C2} and Proposition \ref{P2}(1)).
To solve this equation, we note that, as $(a,1-ax)\in Um(R^{2})$, $a$
divides $1-(1-a-bs)w-by$. Thus, as $R$ is an integral domain, Equation (\ref{EQ8}) has a solution in 
$R^{4}$ iff the system of two equations 
\begin{equation}
1-(1-a-bs)w-by=ta\;\;{and}\;\;t(1-ax)=y[(1-a-bs)z+bx] \label{EQ8}
\end{equation}in the indeterminates $t,x$, $y$, $z$ and $w$ has a solution in $R^{5}$
(this holds even if $a=0$!).

The first equation of System (\ref{EQ8}), rewritten as $w=1+a(w-t)+b(sw-y)$, has the general solution given by $w=1-aq-br$, $t=q+w=1+q-aq-br$ and $y=r+sw=r+s(1-aq-br)$, where $(q,r)\in R^2$ can be arbitrary.

The second equation of System (\ref{EQ8}), rewritten as $t=(at+by)x+(1-a-bs)yz$, has a solution iff $t\in R(at+by)+R(1-a-bs)y$. As $1=(1-a-bs)w+(at+by)$, we have $t\in R(at+by)+R(1-a-bs)y$ iff $t\in R(at+by)+Ry=Rat+Ry$.

From the last two paragraphs it follows that System (\ref{EQ8}) has a solution in $R^{5}$ iff there exists $(q,r)\in R^2$ as mentioned in statement (2). Thus $(1)\Leftrightarrow (2)$.

If we have $(1-a,b)\in Um(R^2)$, then we can choose $(q,r)\in Um(R^2)$ such that $t=1+q(1-a)+rb=0$, hence $t=0\in Rat+Ry=Ry$, and for $e:=1$ and $f:=0$, we have $(a,e),(be+af,1-bs-a)\in Um(R^2)$. Thus to prove that $(2)\Leftrightarrow (3)$ we can assume that $(1-a,b)\notin Um(R^2)$, and hence always $t\neq 0$. 

For $t_1:=\gcd (r-qs,t)$ we write $t=t_1t_{2}$ and $r-qs=\alpha t_1$ with $(\alpha ,t_{2})\in Um(R^{2})$; thus $y=r+s(t-q)=st+\alpha t_1$. We would like to find $(q,r)\in R^{2}$ such that there exists $(g,h)\in R^{2}$ with the property that $gy+aht=t$, i.e., $gst+g\alpha t_1+aht=t$, and thus we must have $g\in Rt_{2}$. Writing $g=t_{2}\beta $ with $\beta \in R$, we want to find $(q,r)\in R^{2}$ such that there exists $(\beta,h)\in R^{2}$ satisfying $1=\beta (st_{2}+\alpha)+ah$, equivalently, such that $(a,st_{2}+\alpha )\in Um(R^{2})$. Replacing $r=qs+\alpha t_1$ in $t$, it follows that $t=1+q-bqs-b\alpha t_1-aq$, hence there exists $\gamma \in R$ such that $1+q(1-bs-a)=t_1\gamma $; we have $t_{2}=\gamma -b\alpha$ and thus $e:=st_{2}+\alpha =s\gamma +(1-bs)\alpha $. As $(\alpha ,t_{2})\in Um(R^{2})$, it follows that $(\alpha ,\gamma )\in Um(R^{2})$. There exists a unique $f\in R$ such that $\alpha =e-fs$ and $\gamma =be+f(1-bs)$; in fact we have $f=t_{2}=\gamma -b\alpha$ and thus $R=R\alpha +R\gamma =Re+Rf$, i.e., $(e,f)\in Um(R^{2})$. There exists $q\in R$ such that $be+f(1-bs)$ divides $1+q(1-bs-a)$ iff $\bigl(be+f(1-bs),1-bs-a\bigr)\in Um(R^{2})$ and hence iff $(be+af,1-bs-a)\in Um(R^{2})$. Thus $(2)\Leftrightarrow (3)$. So Criterion \ref{CR3} holds.

\bigskip \noindent \textbf{Acknowledgement.} The second author acknowledges the support of the research project ``Romanian Hub for Artificial Intelligence - HRIA", Smart Growth, Digitization and Financial Instruments Program, 2021--2027, MySMIS no.\ 334906. The third author would like to
thank SUNY Binghamton for good working conditions. Competing interests: the authors declare none.

\hbox{} \hbox{Grigore C\u{a}lug\u{a}reanu\;\;\;E-mail: calu@math.ubbcluj.ro}
\hbox{Address: Department of Mathematics, Babe\c{s}-Bolyai
University,} 
\hbox{1 Mihail Kogălniceanu Street, Cluj-Napoca 400084, Romania.}

\hbox{} \hbox{Horia F.\ Pop\;\;\;E-mail: horia.pop@ubbcluj.ro} 
\hbox{Address:
Department of Computer Science, Babe\c{s}-Bolyai
University,} 
\hbox{1 Mihail Kogălniceanu Street, Cluj-Napoca 400084, Romania.}

\hbox{} \hbox{Adrian Vasiu,\;\;\;E-mail: avasiu@binghamton.edu} 
\hbox{Address:
Department of Mathematics and Statistics, Binghamton University,} 
\hbox{P.\ O.\ Box
6000, Binghamton, New York 13902-6000, U.S.A.}

\end{document}